\documentstyle[12pt,twoside]{article}
\newtheorem{theorem}{Theorem}
\newtheorem{proposition}{Proposition}
\newtheorem{lemma}{Lemma}

\newtheorem{remark}{Remark}
\newtheorem{corollary}{Corollary}
\newtheorem{definition}{Definition}

\def\R{I\!\!R}
\def\C{I\!\!\!\!C}
\def\N{I\!\!N}

\def\Z{I\!\!\!\!Z}
\def\vsni{\vskip 0.2cm}
\def\var{{\rm var}}

\def\vt{\tilde v}
\def\Var{{\rm Var}}
\def\ui{[0,1]}

\def\A{{\cal A}}

\def\Ft{{\cal{F}}_{\t}}

\def\be{\begin{equation}}
\def\ee{\end{equation}}
\def\qed{\diamondsuit}
\def\O{\Omega}
\def\o{\omega}
\def\t{\theta}

\def\z{\zeta}
\def\l{\lambda}

\def\fc{\tilde{f}}
\def\limsup{\mathop{\overline{\rm lim}}}

\def\S{\Sigma}
\def\s{\sigma}
\def\sp{\sigma^{\prime}}

\def\p{\varphi}

\def\sp{\sigma^{\prime}}
\def\L{\cal{L}}
\def\Nc{\cal{N}}
\def\Mc{\cal{M}}
\def\M{\cal{M}}
\def\P{\cal{P}}
\def\D{\cal{D}}
\def\Dc{\overline {\cal{D}}}

\def\ress{r_{\rm ess}}
\def\ess{\rm ess}
\def\di{\displaystyle}

\begin{document}

\title{On systems with finite ergodic degree}
\date{}

\author{Stefano Isola \thanks{Dipartimento di Matematica e Informatica dell'Universit\`a
di Camerino and INFM, via Madonna delle Carceri, 62032 Camerino, Italy.
e-mail: $<$stefano.isola@unicam.it$>$.}}
\maketitle

\begin{abstract} 
\noindent
In this paper we study the ergodic theory of a class 
of symbolic
dynamical systems $(\O, T, \mu)$ where $T:\O \to \O$ the left shift transformation on $\O=\prod_0^\infty\{0,1\}$
and $\mu$ is a $\s$-finite $T$-invariant measure having the property that one can find a real number $d$ so that
$\mu(\tau^d)=\infty$ but $\mu(\tau^{d-\epsilon})<\infty$ for all $\epsilon >0$, where $\tau$ is the first
passage time function in the reference state $1$. In particular we shall consider invariant measures $\mu$ arising from a potential $V$
which is uniformly continuous but not of summable variation. 
If $d>0$ then $\mu$ can be normalized to give the unique non-atomic
equilibrium probability measure of $V$ for which we compute the (asymptotically) exact mixing rate, of order $n^{-d}$.
We also establish the weak-Bernoulli property and a polynomial cluster property (decay of correlations) 
for observables of polynomial variation.
If instead $d\leq 0$ then $\mu$ is an infinite measure with scaling rate of order $n^d$.
Moreover, the analytic properties of the
weighted dynamical zeta function and those of
the Fourier transform of correlation
functions are shown to be related to one another
via the spectral properties of an operator-valued
power series which naturally arises from a standard inducing
procedure. A detailed control of the singular behaviour 
of these functions in the vicinity of their non-polar singularity at $z=1$ is achieved through an approximation scheme which uses generating functions of 
a suitable renewal process.

\noindent
In the perspective of differentiable dynamics, these are statements about the unique absolutely continuous invariant measure of
a class of piecewise smooth interval maps with an indifferent fixed point. 
\end{abstract}

\section{Introduction}
It is well known that for a subshifts of finite type and equilibrium measure associated to a H\"older continuous
potential $V$ one has exponential decay of correlations. In fact,
H\"older continuity of the potential
implies that $\var_n V$ decays exponentially fast so that the corresponding Ruelle-Perron-Frobenius transfer operator has a spectral gap
when acting
on the Banach space of H\"older continuous functions  (see \cite{Ba1}, \cite{Bo}, \cite{PP}, \cite{Ru4}). 
Moreover, in this case both
the weighted dynamical zeta function and
the Fourier transform of the correlation function 
extend meromorphically to some complex domain where their poles are 
in correspondence with the isolated eigenvalues of the transfer operator
(see \cite{Ba2}, \cite{Hay}, \cite{Pol1}, \cite{Ru3}). If $\var_n V$ decays at a sub-exponential rate
one does not expect a spectral gap any more and the determination of the rate of mixing becomes a 
challenging problem. In \cite{Pol2} this problem has been adressed for potentials of summable variation,
for which it is known \cite{Wal1} that there is only one equilibrium state. In this paper we study
this and related questions for a class of potentials which are not even of summable variation, but
have an induced version which is H\"older continuous. In this case, there is a $\sigma$-finite 
invariant measure which is either infinite or can be
normalized to give a (non-unique) equilibrium state, depending on the value of a parameter that 
we identify as the {\sl ergodic degree}. Roughly speaking, this parameter controls in a continuous fashion the number of finite moments
possessed by the first passage time function in a given reference set. 

\noindent
The main motivation for this study comes from
an attempt to understand the ergodic properties
of some class of non-uniformly hyperbolic dynamical systems, the simplest example being that of smooth interval maps which are expanding
everywhere but at an indifferent fixed point (see \cite{Th}, \cite{PS}, \cite{LSV1}, \cite{Yo}, 
\cite {Yu}, \cite{MRTVV}, \cite{Hu}). 

\noindent
The paper is organized as follows. In Section \ref{inducing} we give some preliminaries on the inducing procedure and
the main assumptions are settled down. In Section \ref{opvalpowser} we introduce 
an operator-valued power series ${\M}_z$ which will play an essential role in the sequel and study its spectral properties when acting
on a Banach space of locally H\"older continuous functions.
An algebraic relation between ${\M}_z$ and the transfer operator $\L$ of the original system is established in Section
\ref{section4}. 
This relation is then used to construct the $\s$-finite invariant measure mentioned above and to discuss some of
its properties, depending of its ergodic degree $d$. Moreover, the logarithm of the leading 
eigenvalue of ${\M}_z$, for $z$ in a complex open neighbourhood of $(0,1]$,
is interpreted as a pressure function $P(z)$
for a suitable `grand canonical' potential, and in Section \ref{pressione} we examine its behaviour, showing in particular how the number of its derivatives at $z=1$ is related to the ergodic degree $d$. 
An approximation scheme based on a renewal Markov chain is introduced in Section \ref{markovapprox} and the asymptotic
behaviour when $z\uparrow 1$ of the operator valued function $(1-{\M}_z)^{-1}$ acting on suitable test functions is determined in terms of
the generating function
of this renewal process. 
In Section \ref{rensca} we study the Fourier transform of a `renewal density' sequence
yielding the probability to observe a return in the reference set after $n$ iterates, and we show that 
it is asymptotically equivalent to the corresponding renewal sequence for the Markov approximation. Thereafter,
this is extended to any Borel set where the first passage time function in the above reference set is bounded, thus determining
the (asymptotically) exact scaling rate when the measure is infinite ($d\leq 0$) and the 
(asymptotically) exact mixing rate when it is finite ($d>0$). In the latter case, we also establish a weak-Bernoulli
property and a polynomial decay of correlations for test functions of polynomial variation (Section \ref{poldecay}).
In Section \ref{zetafunc} we study dynamical zeta functions. We first establish an algebraic relation between zeta functions
of the original and the induced system which is the counterpart of the operators relation mentioned above. We then
show that the singular behaviour of the dynamical zeta function is characterized by a non-polar singularity at $z=1$ 
which can be inspected in terms of
the generating function of the approximating renewal process.
Finally, in Section \ref{app} we apply the preceeding results to the symbolic description of the dynamics of interval maps with
indifferent fixed points. In this simple situation we have a nice control of the asymptotic behaviour of the first passage function
which allows to partially sharpen the general results. Although we believe that the approach described here can be carried out for
more general examples of non-uniformly hyperbolic systems, we leave the actual investigation of this point to be discussed elsewhere.
\vsni  
\noindent
\section{Inducing on shift spaces}\label{inducing}
We shall consider the following situation.
Let $\O=\prod_0^\infty\{0,1\}$ be 
the set of all one-sided sequences
$\o = (\o_0\o_1\dots )$, $\o_i\in\{0,1\}$ and denote the left shift map on $\O$ by $T$.
We shall take the state $1$ as a reference set and consider the 
{\sl first passage function} $\tau : \O \to \N$ given by:
\be\label{fpf}
\tau(\o)=1+\inf\{i\geq 0 \, :Ê\, \o_i=1\},
\ee
(with $\inf \emptyset = \infty$).
The levelsets
\be
A_k:=\{\o \in \O : \tau (\o)=k\},\qquad k\geq 1,
\ee
are both open and closed and will be called {\sl partition sets}, as they form
a partition of $\O_0:=\O\setminus \{0^\infty\}$, where $0^{\infty}$
denotes the singleton $(000\dots )$. 
Let $0^{n}1$ denote the word ${\underbrace{0\dots 01}_{n+1}}$.
The $A_k$'s will be sometimes denoted as cylinder sets:
\be\label{cylinder}
A_k\equiv [0^{k-1}1]:=\{\o \in \O \; | \; \o_i=0, \,0\leq i <k-1,\; \o_{k-1} =1\}.
\ee
\noindent
Let moreover
\be \tau_0 \equiv \tau,\quad\hbox{and}\quad \tau_k=
\inf\{i>\tau_{k-1}\, :Ê\, \o_i=1\},\quad k>0,
\ee
be
the
{\it sequence of successive entrance times in} 1, and let
\be
\label{passages0}
\s_0\equiv \tau,\quad\hbox{and}\quad \s_k=\tau_{k}-\tau_{k-1},\quad k>0,
\ee
be the {\it sequence of times between passages}. Both sequences are infinite on the
residual subset $\O\setminus \cup_{k=0}^\infty T^{-k}0^\infty$. 
Now, given a function $U: \O \to \C$, we define two different notions of variation:
\be
\var_n U  =  \sup \biggl\{ \, | U(\o)-U(\o') |\,:\, 
\inf \{k\geq 0\,:\, \o_k\ne \o'_k\}=n \biggr\},
\ee
and
\be
 \Var_n U  =  \sup \biggl\{ \, | U(\o)-U(\o') |\,:\, 
\inf \,\{k\geq 0\,:\, \s_k(\o)\ne \s_k(\o')\}=n \biggr\}.
\ee
We then say that $U$ is {\sl uniformly continuous} if $\forall \, \epsilon >0$
there exists $n\geq 1$ such that $\var_n U  < \epsilon$. On the other hand
$U$
will be called {\sl locally H\"older continuous}
if there is a constant $C>0$ 
such that $\forall n \geq 1$, $\Var_n U \leq  C\, \t^n$ (notice that nothing is required for
$\Var_0\, U$).
Let $V:\O \to \C$ be given
and define its {\sl induced version} $W:\O\to \C$ as
\be
\label{inducedversion}
W(\o) = \sum_{k=0}^{\tau(\o)-1}V(T^k\o).
\ee
Conversely, we have
\be \label{converse0}
V(\o) = \cases{W(\o) &if $\o\in A_1$,\cr
               W(\o)-W(T\o) &if $\o\in A_k, \; k>1$,\cr}
\ee
We shall be interested in examples of uniformly continuous maps $V:\O \to \R$
whose induced version $W:\O\to \R$ is locally H\"older continuous. 
Whenever the function $V$ is viewed as a {\sl potential} corresponding
to the map $T:\O\to \O$, its induced version $W$ is naturally interpreted as
the potential for the {\sl induced map} $T^\tau : \O \to \O$.

\noindent
We shall now introduce a class 
of real valued maps $V$  on the space $\O$ whose regularity properties are dictated by
those of their induced version $W$.

\vskip 0.5cm
\noindent
{\bf PROPERTIES.} {\sl 
We will be considering potential functions $V:\O \to \R$ having the following properties:}
\begin{enumerate}
\item  {\sl $V(\o)$ is continuous and satisifes $-C_1 \leq V(\o)\leq 0$ (for 
some positive constant $C_1$), with
$V(0^{\infty})=0$; }
\item {\sl  there is a constant $C_2 >0$ s.t.}
$$
\sum_n \; \exp \left( {\sup_{\o \in A_n}W(\o) } \right) < C_2;
$$
\item
{\sl there is $0<\t <1$ and a constant $C_3>0$ such that for any $n \geq 1$,}
$$
\Var_n W \leq C_3\, \t^{n};
$$
\item {\sl the limit
$$
{\wp}(V):=\lim_{n\to \infty} {1\over n}
\log \, \sum_{\o_0,\o_1,\dots ,\o_{n-1}}\,\,\exp \left(\sup_{\o'}
{\sum_{i=0}^{n-1} V(\o_i\dots \o_{n-1}\o')}\right)
$$
exists and is equal to $0$.}
\end{enumerate}
\begin{remark}
{\rm One readily obtains from (\ref{converse0}) that
\be
\var_n V = \sup_{\o \in A_{n+1}}| V(\o)|  ,
\ee
so that Property 1 entails that $V$ is uniformly continuous and, in particular, $\var_n V$ eventually decreases
monotonically. On the other hand, by Property 2, $V$ is not of summable variation:
$$
\sum_{k< n}\var_{k} V \geq \inf_{\o \in A_n} |W(\o)| \to \infty\quad\hbox{as}\quad n\to \infty.
$$
Moreover, by Property 3, its induced version $W:\O \to \R$ is 
locally H\"older continuous, even though not bounded from below. }
\end{remark}
\begin{remark}
{\rm The function $\wp (V)$ defined in Property 4 is called the (topological)
pressure of $V$. Its existence follows from Property 1 along with the sub-additivity of the 
sequence $\{\, \sup_{\o'}{\sum_{i=0}^{n-1} V(\o_i\dots \o_{n-1}\o')} \, \}_{n\in \N}$. Notice that
for any real number $c$ we have $\wp (V+c)=\wp (V)+c$, and therefore $\wp (V-\wp (V))=0$.
Hence, it always possible to reduce to the case of potentials with zero pressure.}
\end{remark}
{\bf Example 1.} Let $1=p_0>p_1>p_2>\cdots$ be a sequence of real 
numbers such that $p_k>0$, $\sum_{k\geq 1} p_k =1$ and $q_k =\log (p_{k-1}/p_k)\searrow 0$. Define a continuous function 
$V:\O \to \R$
as
$$
V(\o) =\cases{-q_k &if $\o\in A_k$,\cr
              \;\;\; 0                  &if $\o=0^{\infty}$,\cr
              \;\;\; 0                  &otherwise.\cr}
$$
One easily checks that $\var_n V = q_{n+1}$,
so that $V$ is uniformly continuous on $\O$. Notice that $V$ is not of summable variation:
$\sum_{k\leq n} q_{k} =-\log \, p_n \to \infty$ as $n\to \infty$.
On the other hand $W(\o)=-\sum_{k=1}^{\s_0}q_k=\log p_{\s_0}$ so that $\Var_n W = 0$, $\forall n \geq 1$. 
Notice that using the correspondence $A_k \leftrightarrow k$ we can 
relate this example to a Markov chain with state space $\N$
and transition probabilities 
\be 
p_{ij}=\cases{  p_j, &if $\; \; i=1, j\geq 1$, \cr
               1, &if $\; \; 1\leq j =i-1$,\cr 
              0,  &otherwise.}
\ee
To see this, we set $x_i=k$ if $T^i\o \in A_k$ and 
$x=(x_0x_1\cdots)\in X=\N^{\N_0}$, and note that there is a one-to-one
correspondence
between points $\o \in \O_0$ and points $x\in X$ satisfying the
compatibility condition: given $x_i$ then either $x_{i-1}=x_i+1$ or $x_{i-1}=1$.
Such a correspondence induces a natural action of the shift $T$ on $X$. 
With slight abuse of notation, we can thus write
the function $V$ as a function of $x$, which turns out to depend
only on the first two coordinates:
$$
V(x) =\log \left({p_{x_0}\, p_{x_0x_1}\over p_{x_1}}\right).
$$
This can be viewed as a `two-body potential' function for the Markov shift $(X,T)$. 
In particular, there is no mutual interaction between the `spins' $\s_k$
given by the times between passages in the state $1$: the sequence $\{\s_k\}$
is isomorphic to a sequence of i.i.d.r.v. and 
$\tau_n =\s_0+\cdots +\s_n$ is a stationary  renewal process
under the probability measure $P(\s_k=\ell)=p_{\ell}$.

\section{An operator-valued power series} \label{opvalpowser}
We define the operation of {\it inducing} as the bijection $\iota : \O_0 \to \S$
where $\S := \iota (\O_0)$ is the non-compact set of all sequences 
$\s = \iota (\o)$ given by times between passages, i.e. $\iota (\o) = \s_0(\o)\s_1(\o)\dots$.
The map $S:=\iota^*\, T^\tau$ acts as
a left shift on $\S$. 
In this Section we shall exclusively work with objects (functions, measures, operators) living on the
symbol space $\S$, temporarily forgetting its origin as an induced space from $\O_0$ through the
map $\iota$.
We denote by $\Ft(\S)$ the Banach space of complex valued functions on $\S$
which are finite with respect to the norm
$\Vert U \Vert_{\t} = |U |_{\infty}+|U |_{\t}$, where
$|U |_{\t}$ is the least Lipschitz constant of $U$ wrt the metric 
\be
d_\theta (\s,\s') = \t^{\di\, \inf \,\{j\geq 0\,:\, \s_j\ne \s'_j\}}.
\ee
that is 
\be
|U |_{\t}= \sup \left\{  {|U (\s)-U (\s')|\over \t^n}\,:\, 
 n=\min \{j\geq 1\,:\, \s_j\ne \s'_j\}\right\}.
\ee
We shall denote with the same symbol $W$ the projection $\iota^*W:\S\to \R$ of the induced
potential defined in (\ref{inducedversion}).
Notice that even though $W$ is not bounded from below,
the function 
\be
\label{Boltzmannfactor}
{\psi}(\s):= \exp W (\s)
\ee
clearly is. In particular, by Property 2 we have that 
\be
\label{stima1}
\sum_{S \s'=\s}\psi (\s')=\sum_{k=1}^{\infty} \psi (k\s) < C_2,
\quad\hbox{for all}\quad \s\in \S.
\ee
Here $k\s$ denotes the sequence $(k\s_0\s_1\dots )\in \S$. 
This and Property 3 entail that $\psi \in \Ft(\S)$. 
Let $z\in \C$ and ${\M}_{z} : {\Ft} (\S) \to {\Ft} (\S)$ be
the operator-valued power series defined as follows (see \cite{PS} for related objects):
\be
\label{operator}
{\M}_{z} = \sum_{n=1}^{\infty}\, z^n{\Mc}^{(n)},\qquad
({\Mc}^{(n)}\, f)(\s) = \psi(n\s)\,  f(n\s).
\ee
Alternatively,
one can think of ${\M}_{z}$ as the transfer operator associated
to the `grand canonical' complex potential
\be
W_z(\s) := W(\s) + \s_0 \log z
\ee
where we take the determination of $\log z$
which is real for $z>0$.
\vsni
\noindent
\begin{lemma} \label{radius} The power series of 
${\M}_{z}$ when acting on ${\Ft}$ has radius of convergence
bounded from below by $1$ and, moreover, 
it converges absolutely at every point of the unit circle.
\end{lemma}
{\it Proof.}
The radius of convergence of ${\M}_{z}$ is  
$\lim_{n\to \infty}\Vert {\Mc}^{(n)} \Vert_{\t}^{-1/n}$.
We have
$$
|{\Mc}^{(n)} f(\s)| \leq \psi(n\s) |f|_{\infty}
$$
and also
$$
|{\Mc}^{(n)} f(\s) - {\Mc}^{(n)} f(\s') | \leq \,
\psi(n\s)\, \left( \, \t \, |f|_{\t} + C \, \t \, |f|_{\infty}\, \right).
$$
Hence, 
\be\label{estim}
\Vert {\Mc}^{(n)} \Vert_{\t} \leq C \, \sup_{\s}\, \psi(n\s)
\ee
and the assertion follows from (\ref{stima1}). $\qed$
\vskip 0.1cm 
\noindent
For any fixed  $z\in \ui$ we now set
\be
\label{Lambda_n}
\Lambda_n(z) := \sup_{\s \in \S} \sum_{k_1=1}^{\infty}\dots \sum_{k_n =
1}^{\infty} z^{k_1+\dots + k_n} \prod_{i=1}^n\psi(k_i\dots k_n\s)
\ee
and
\be
\label{pressure}
P(z) := \lim_{n\to \infty}{1\over n}\log \Lambda_n (z).
\ee
One easily checks that under our hypotheses on the the potential $V$ 
(and for each $z$ as above) the
sequence $(\Lambda_n(z))_{n\in \N}$ is sub-multiplicative, and
therefore the limit (\ref{pressure}) exists and satifies $-\infty \leq P(z) < \infty$.
Also notice that, using
(\ref{inducedversion}) and  (\ref{Boltzmannfactor}) with $\s = \iota (\o)$, we can write
$$
\prod_{i=1}^n\psi(k_i\dots k_n\s) = 
\prod_{j=0}^{k_1+\cdots +k_n-1}\p(T^j\,{\underbrace{0\dots 01}_{k_1}}\;\cdots \;
{\underbrace{0\dots 01}_{k_n}}\,\, \o), 
$$
where $\p = \exp V$. Property 4 then implies that
$$
P(1)=\lim_{n\to \infty} {1\over n} \log \Lambda_n(1) = 0.
$$
We shall however say more. In the next theorem we prove that for any fixed  $z\in \ui$
the function $\exp P(z)$ is equal to the spectral radius $r({\M}_{z})$ of ${\M}_{z}:{\Ft}\to {\Ft}$.
The monotonicity of $z \to {\M}_z$ for $0\leq z \leq 1$ thus implies that 
$\exp P(z)$ is strictly increasing, ranging from $0$ to $1$ when $z$ ranges from $0$ to $1$.
\begin{remark} \label{markovpressure}
{\rm For $V$ and $W$ as in Example 1
we find $\Lambda_n(z)=\left(\sum_{k=1}^{\infty}z^k\,p_k\right)^n$ so that 
$\exp P(z)$ is but the generating function of the numbers $p_k$.}
\end{remark}
For $z$ in some complex neighbourhood of $(0,1]$
the quantity $P(z)$ will be interpreted as
the {\sl pressure} associated to the potential $W_z$. 
Before stating the next result we let
$\D$ and $\Dc$ denote the open unit disk $\{z\, :\, |z|< 1\}$ and 
its closure $\{z\, :\, |z|\leq 1\}$, respectively.
Moreover, we
recall that 
the spectrum of a bounded linear
operator $K$ can be decomposed into a discrete part, made up of isolated
eigenvalues of finite multiplicity, and its complement, the essential
spectrum, denoted by $\ess (K)$. 
The essential spectral radius
is then defined as $\ress (K) = \sup \, \{ \, |\lambda| \, : \, \lambda
\in  \ess (K) \, \}$. 
\begin{theorem} \label{spectrum}
Let $z\in \Dc$ and
${\M}_{z}$ be acting on ${\Ft}$.
\begin{enumerate}
\item The spectral radius $r({\M}_{z})$ is bounded above by
$\exp P(|z|)$.
\item The essential spectral radius $\ress ({\M}_{z})$ is
bounded above by $\t \, \exp P(|z|)$.
\item There is at most one eigenvalue of modulus $\exp P(|z|)$,
which is simple, and exactly one at $\exp P(z)$ if 
$z$ is real and positive. The rest of the spectrum is contained in a disc
of radius strictly smaller than $\exp P(|z|)$. In addition, if $|z|\leq 1$ but
$z\ne 1$ then $1$ is not an eigenvalue of ${\M}_{z}$.
\end{enumerate} 
\end{theorem}
{\it Proof.} 
It is easy to check that for all $f\in {\Ft}$
\be
\label{stimabase1}
|{\M}_{z}^n f|_{\infty}\leq \Lambda_n(|z|)\, |f|_{\infty} \leq \Lambda_n(|z|)\,
\Vert f \Vert_{\t}
\ee
and 
\be
\label{stimabase2}
|{\M}_{z}^n f|_{\t}\leq \Lambda_n(|z|)\, (C\, |f|_{\infty}+\t^n |f|_{\t})
\leq \Lambda_n(|z|)\, (C +1)\Vert f \Vert_{\t}.
\ee 
Therefore $\Vert {\M}_{z}^n \Vert_{\t} \leq (C+2)\Lambda_n(|z|)$, where
we have also denoted by $\Vert \,\,\, \Vert_{\t}$ the operator norm.
Thus, the spectral radius formula implies that 
$$
r({\M}_{z})\leq \exp P(|z|).
$$
To estimate the essential spectral radius 
we now extend an argument given in \cite{Ke} (see also \cite{Pol1}) to the
present infinite-alphabet situation.
For any $n\geq 1$, let $\S^{(n)}$ be the set of all words $\eta$
of length $n$, i.e. words of the form 
$\eta = (\s_0 \dots \s_{n-1})$. Moreover, given $N>0$, set 
$G_n^N=\{\eta\in \S^{(n)} \, : \, \s_i < N, \, 
0\leq i<n\, \}$ and $B_n^N=\S^{(n)} \setminus G_n^N$. Now choose
$N=N(n)$ in such a way that
\be
\label{stima2}
\sup_{\s\in \S}\,\sum_{\eta \in B_n^N}\, |z|^{S_n(\eta\s)}\,e^{S_nW (\eta\s)} 
\leq C\, \Lambda_n(|z|)\theta^n
\ee
where $S_n(\s)= \sum_{j=0}^{n-1}\s_j$ and $S_nW = \sum_{j=0}^{n-1}W\circ S^j$.
Having fixed $n$ and $N$,
let $E_N^n$ be the finite rank operator acting on $f\in \Ft$ as
follows:
$$
(E_N^n f)(\s) = \sum_{\eta \in  G_n^N}f_\eta\, \chi_{\eta}(\s)
$$
where $\chi_{\eta}(\s)$ is the characteristic function of the cylinder set
$D_\eta =\{\s \in \S \, :Ê\, \s_j=\eta_j, \, j=0,\dots ,n-1\}$ and $f_\eta$
is the value of $f$ at some (arbitrarily) chosen point of $D_\eta$.
Then, for any $0<|z|\leq 1$ and for any pair $n,N>0$, the operator 
$K_{z,N}^n={\M}^n_{W_z}E_N^n$ is compact. 
Put moreover 
$$
\fc(\s)=\cases{(f-E_N^n f)(\s),&if $[\s_0, \dots ,\s_{n-1}] \in G_n^N$; \cr
f(\s),&otherwise. 
\cr}
$$
It is easy to see that 
$$
\sup_{\eta \in G_n^N} |\fc(\eta\s)|\leq C\,|f|_{\t}\, \t^n
\quad\hbox{and}\quad 
\sup_{\eta \in G_n^N} |\fc(\eta\s)-\fc(\eta\sp)|\leq C\,|f|_{\t}\, \t^n
$$
Hence,
\begin{eqnarray}
& & |\, ({\M}^n_{z} - K_{z,N}^n)f(\s)\, | = 
|\, {\M}^n_{z}\fc(\s)\, | \leq \nonumber \\
& &\leq | f|_{\infty}
\sup_{\s\in \S}\,\sum_{\eta \in B_n^N}\, |z|^{S_n(\eta\s)}\,e^{S_nW(\eta\s)} 
+\Lambda_n(|z|)\, \sup_{\eta \in G_n^N} |\fc(\eta\s)|
 \nonumber  
\end{eqnarray}
Therefore, using (\ref{stima2}), we get
$$
|\, ({\M}^n_{z} - K_{z,N}^n)f\, |_{\infty}\leq  C\, \Vert f\Vert_{\t}\,
\t^n\, \Lambda_n(|z|).
$$
We now estimate the variation. Let $\s, \sp \in \S$ be two sequences having 
the first $k$ symbols in common. Then we have
\begin{eqnarray}
& &|\, {\M}^n_{z}\fc(\s)-{\M}^n_{z}\fc(\sp)\, | \leq \nonumber \\
& &\leq 
\sum_{\eta \in G_n^N}\, 
|z|^{S_n(\eta\s)}\,\left( e^{S_nW(\eta\s)}\fc(\eta\s)-
e^{S_nW(\eta\sp)}\fc(\eta\sp)\right) +
\nonumber \\
& &+\sum_{\eta \in B_n^N}\, 
|z|^{S_n(\eta\s)}\,\left( e^{S_nW(\eta\s)}f(\eta\s)-
e^{S_nW(\eta\sp)}f(\eta\sp)\right)
={\rm I}+{\rm II}. \nonumber
\end{eqnarray}
Now a routine calculation shows that 
$$
{\rm I} \leq C \, \Lambda_n(|z|) \, |f|_{\t}\, \t^{n+k}.
$$
Moreover
\begin{eqnarray}
& &{\rm II} \leq \sum_{\eta \in B_n^N}\, 
|z|^{S_n(\eta\s)}\,e^{S_nW(\eta\s)}\,\left( f(\eta\s)-
f(\eta\sp)\right)+ \nonumber \\ & &+
\sum_{\eta \in B_n^N}\, 
|z|^{S(\eta\s)}\,e^{S_nW(\eta\s)}\, f(\eta\sp)\, \left(1-
e^{S_nW(\eta\sp)-S_nW(\eta\s)}\,\right) \nonumber \\
& &\leq C\, \Lambda_n(|z|) \, |f|_{\t}\, \t^{n+k}
+C\, \t^k \,  |f|_{\infty}\,
\sum_{\eta \in B_n^N}\, |z|^{S_n(\eta\s)}\,e^{S_nW(\eta\s)} \nonumber 
\end{eqnarray}
so that, using (\ref{stima2}), we get 
$$
|\, ({\M}^n_{z} - K_{z,N}^n)f\, |_{\t}\leq  C\, \Vert f\Vert_{\t}\,
\t^n\, \Lambda_n(|z|).
$$
Putting together the above and Nussbaum formula \cite{Nu}
we have thus proved that 
$$
\ress ({\M}_{z}) \leq \t \, \exp P(|z|).
$$
We are now going to prove the third statement. 
Let us first notice that, if $z$ is real and positive, we have
$$
r({\M}_{z}) =\lim_{n\to \infty} 
\left( \, \Vert {\M}^n_{z} \Vert_{\t}\, \right)^{1/n}
\geq \lim_{n\to \infty} \left(\, |{\M}^n_{z}1 |_{\infty}\, \right)^{1/n} = 
\exp P(z)
$$
and therefore $r({\M}_{z}) = \exp P(z)$. 
Next, let us see that for $z$ real and positive
$\exp P(z)$ is a maximal simple eigenvalue of ${\M}_{z}$
and the remainder of the spectrum is contained in a disk of radius 
strictly smaller than $\exp P(z)$. 
To this end, we proceed as in (\cite{CI}, Theorem 2.1) and construct a
sequence of compact spaces $\S_N$, $N\in \N$, whose elements are 
sequences $\s=(\s_0\s_1\dots )$ with $\s_j\in\{1,\dots, N\}$.
Clearly, $\S_N \subset \S_{N+1} \subset \dots \subset \S$. 
For any $0<z\leq 1$ define a family of operators ${\M}_{z,N}:\Ft(\S_N)\to
\Ft(\S_N)$ by
$$
{\M}_{z,N}= \sum_{k=1}^{N} z^k\, {\Mc}^{(k)}.
$$
Now, having fixed $N\in \N$ and $z \in (0,1]$, a 
Ruelle-Perron-Frobenius theorem holds for ${\M}_{z,N}$ 
(see, e.g., \cite{PP}, Theorem 2.2).
Let $\lambda_{z,N}$ be the simple eigenvalue with largest modulus
of ${\M}_{z,N}$ and 
${h}_{z,N}\in {\Ft}(\S_N)$, $\nu_{z,N}\in {\Ft}^*(\S_N)$ be such that 
$$
{\M}_{z,N}{h}_{z,N}=\lambda_{z,N}{h}_{z,N}\quad\hbox{and}\quad
{\M}^*_{z,N}{\nu}_{z,N}=\lambda_{z,N}{\nu}_{z,N}.
$$
Putting 
$$
\Lambda_{n,N}(z) = \sup_{\s \in \S_N} \sum_{k_1=1}^{N}\dots \sum_{k_n =
1}^{N} z^{k_1+\dots + k_n} \exp \sum_{i=1}^nW (k_i\dots k_n\s)
$$
and 
\be
\label{P_N}
P_N(z) = \lim_{n\to \infty}{1\over n}\log \Lambda_{n,N}(z)
\ee
we find 
\begin{eqnarray}
\log \l_{z,N} &=& \limsup_{n\to \infty}
{1\over n}\log {\M}_{z,N}^n{h}_{z,N}  \nonumber \\
&=&\limsup_{n\to \infty}
{1\over n}\log {\M}_{z,N}^n 1 =P_N(z) \nonumber 
\end{eqnarray}
because $D_{z,N}^{-1} \leq {h}_{z,N} \leq D_{z,N}$ for some
positive constant $D_{z,N}$.
Now, reasoning as in 
(\cite{CI}, Theorem 2.1) one then shows that, having fixed $0<z\leq 1$, the triple
$(\, \lambda_{z,N}, {h}_{z,N}, {\nu}_{z,N}\, )$ converges uniformly
to the triple $(\, \lambda_z, {h_z}, \nu_z\, )$ which is uniquely 
determined by the conditions
\begin{eqnarray}\label{eigen}
&&\lambda_z >0, \quad h_z>0, \quad\int h_z\, d\nu_z=1\\
&&{\M}_{z}{h_z}=\lambda_z{h_z},\quad
{\M}^*_{z}{\nu_z}=\lambda_z{\nu_z}. \nonumber
\end{eqnarray}
Clearly we have $\lambda_z= \exp P(z)$ where $P(z)$ is defined in (\ref{pressure}). 
Furthermore, 
\be
\label{convergence}
\lambda_z^{-n}{\M}^n_{z} f\to h_z \int f\, d\nu_z
\quad\hbox{uniformly}\quad\forall \, f\in \Ft (\S).
\ee
The fact that the rest of the spectrum
is contained in a disk of radius strictly smaller than $\exp P(z)$
now follows from the argument given in 
(\cite{PP}, p.26), which relies on the basic 
inequalities (\ref{stimabase1}) and (\ref{stimabase2}). 

\noindent
Finally, in order to deal with complex values of $z$, let us write 
$z=|z|e^{i\phi}$ and $W_z= a + ib$ where
$a(\s) = W(\s) + \s_0 \log |z|$ and $b(\s) = \s_0\phi$. Proceeding as in
(\cite{PP}, Chapter 4), one then shows that if 
${\M}_{z}$ has an eigenvalue of modulus
$\exp P(|z|)$, then 
${\M}_{z}= \vartheta \, {\cal S} {\M}_{a}{\cal S}^{-1}$
where ${\M}_{a}={\M}_{{|z|}}$, ${\cal S}$ is a multiplication operator 
and $\vartheta \in \C$, 
$|\vartheta |=1$, so that the spectral properties of ${\M}_{z}$
follow from those of ${\M}_{a}$ and hence from the above discussion; 
if, instead,
${\M}_{z}$ has no eigenvalues of modulus
$\exp P(|z|)$, then  its spectral radius is strictly smaller than $\exp P(|z|)$ and $\exp{(-nP(|z|))}\,{\M}^n_{z} \to 0$ in the 
$\Vert \,\, \Vert_{\t}$-operator topology.  To complete the proof 
of the Theorem, suppose
there exist $z=e^{i\phi}$ with $\phi ({\rm mod}2\pi) \neq 0$ 
and $h_z\in {\Ft}(\S,\C)$ such that
${\M}_{z}h_z=h_z$ (the case $|z|<1$ follows by just observing that $\exp P(z)$ is strictly
increasing in $z$ when $0\leq z \leq 1$ and $r({\M}_{z})\leq \exp P(|z|)$).
Let ${\cal V}_z: {\Ft}(\S,\C) \to {\Ft}(\S,\C)$ be the operator defined by the relation
\be
{\M}_z (f\cdot h_1) = h_1 \cdot {\cal V}_z f \, .
\ee
Note that ${\cal V}_11=1$. Clearly ${\M}_z$ and ${\cal V}_z$ have the same spectrum and
all eigenvalues have the same geometric multiplicities.
Under the above supposition for $z=e^{i\phi}$ we would have ${\cal V}_zg_z =g_z$ with $g=h_z/h_1$, and
from ${\cal V}_11=1$ and ${\cal V}_1(e^{i\,\phi \, \s_0} \, g_z)=g_z$ we would obtain
${\cal V}_1(e^{i\,\phi \, \s_0} \, g-g\circ S)=0$. 
Iterating this argument we would get a function $g_z\in {\Ft}(\S,\C)$ such that
$g_z(S^n \s)=g_z(\s)e^{i\,\phi \,\sum_{k=0}^{n-1} \s_k} + r(\s)$
with $r(\s)\in {\rm ker}\, ({\M})$. This would give
\be
\int_\S g_z(S^n\s) \, d\nu_1 (\s) = \int_\S g_z(\s)\, e^{i\,\phi \,\sum_{k=0}^{n-1} \s_k} \, d\nu_1 (\s) .
\ee
However $S:\S \to \S$ is mixing so that 
we have found a contradiction. 
$\qed$

\vsni
\noindent
Some interesting consequences for the ergodic theory of the shift
map $S$ on $\S$ (together with the weight function $\psi$),
can be obtained by putting $z=1$ in the above theorem. Related results can be found in 
\cite{LSV2}, \cite{MD}, \cite{Br}, \cite{Sa1}.

\noindent
First, we shall say that a measure $\nu$ is {\sl conformal} if it
satisfies
\be
\label{conformal}
\int_\S {\M}f\, d\nu = \int_\S f\, d\nu,\qquad \forall f\in \Ft.
\ee
Now, from the above theorem it follows that the
bounded linear operator ${\M}:\Ft \to \Ft$ given by ${\M}\equiv {\M}_{1}$, 
that is the usual transfer operator
associated to the shift map $S$ on $\S$
and the weight $\psi$, has spectral radius $r({\M})=1$. Let $\nu \equiv \nu_1$, $h \equiv h_1$
$\lambda \equiv \lambda_1$ with $\nu_1$, $h_1$, $\lambda_1$ defined in eq (\ref{eigen}).
Then the element $\nu$ of the dual space $\Ft^*(\S)$
is a conformal measure. 
Furthermore, we recall that a probability measure $\rho$ on $\S$
is called a Gibbs measure (in the sense of Bowen, see \cite{Bo}) 
if there exists $\Psi \in {\cal C}(\S)$ such that
$$
A\leq {\rho ([\s_0\cdots \s_{n-1}]) \over e^{-n\, C+ \sum_{k=0}^{n-1}\Psi (S^k\s)} }\leq B
$$
for $n>0$ and fixed constants $A,B>0$ and $C\in \R$. Here
$[\s_0\cdots \s_{n-1}]$ denotes the cylinder set 
$\{\s' \in \S \, :Ê\, \s'_i=\s_i,\, 0\leq i < n\}$ and $\s$ is any point in it. 
\begin{corollary} \label{gibbs} The measure $\rho :=h \cdot \nu$ is a $S$-invariant Gibbs measure for
$C=0$ and $\Psi = W$. Moreover it is uniformly mixing, i.e.  there is a constant $0<\vartheta <1$ 
such that for any pair of cylinders
$E=\{\s \in \S : \s_i=e_i, 0\leq i \leq s\}$ and $F=\{\s \in \S : \s_i=f_i, 0\leq i \leq r\}$ 
we can find a constant $M=M(E,F,\vartheta)$ such that 
$$
\left|{\rho(S^{-n}F \cap E)\over \rho (E)} - \rho (F)\right|\leq M\, \vartheta^n.
$$
\end{corollary}
{\sl Proof.} We have
\begin{eqnarray}
\rho ([\s_0\cdots \s_{n-1}]) &=& \int_\S \chi_{[\s_0\cdots \s_{n-1}]}(\s)\cdot h (\s) d\nu (\s)
\nonumber \\
&=&\int_\S {\M}^n (\chi_{[\s_0\cdots \s_{n-1}]} \cdot h )(\s) \, d\nu (\s) \nonumber \\
&=&\int_\S  h(\s_0\cdots \s_{n-1}\s) \psi (\s_0\cdots \s_{n-1}\s)\, d\nu (\s) \nonumber
\end{eqnarray}
Now, since $W$ is locally H\"older continuous,
$\sum_{n\geq 1}\Var_n W \leq C_3\t/(1-\t)< \infty$. Therefore, for any pair $\s,\s'\in \S$,
we have 
\be
e^{{-C_3\theta\over 1-\t}}\, {\M}^n1(\s') \leq {\M}^n1(\s) \leq e^{{C_3\theta\over 1-\t}}
{\M}^n1(\s') \, .
\ee
Taking limits
as $n\to \infty$ and recalling that ${\M}^n 1\to h$ we get $\Vert \log h\Vert_\infty < \infty$ so that we can find a constant $D>0$ s.t.
$D^{-1} < h < D$. The assertion, with
$B=A^{-1}=D\, e^{C_3\t/(1-\t)}$, now follows from the above identity. The $S$-invariance follows
from the fact that $h$ is a fixed function for $\M$
and the uniform mixing property from the existence of a spectral gap for 
${\M}:\Ft \to \Ft$ (see \cite{Ru4}). $\qed$

\vskip 0.2cm
\noindent
We end this Section by studying the operator-valued function $z\to (1-{\M}_z)^{-1}$. We first
recall that $z\to {\M}_z$ is holomorphic in $\D$ and continuous on $\Dc$.
Therefore, if $(1-{\M}_z)^{-1}$ exists (as a bounded operator acting on $\Ft$) for all
$z$ in some open subset $D\subseteq \D$ then  $z\to (1-{\M}_z)^{-1}$ is holomorphic in $D$.
We know from Theorem \ref{spectrum} that the spectral radius of ${\M}_z:\Ft \to \Ft$ is 
bounded above by $\exp P(|z|)$ which is $<1$ for $0\leq |z| <1$.
In addition, from the last statement of Theorem \ref{spectrum} 
it easily follows that for $|z|\leq 1$ and $z\ne 1$ there is no eigenvalue of modulus $1$ (see also Proposition \ref{relations} below). This shows that
the function $z\to (1-{\M}_z)^{-1}$ is holomorphic in $\D$ and extends continuously to
$\Dc \setminus \{1\}$.

\noindent
Now, standard analytic perturbation theory 
(\cite{Ka}, Section 7.1)
implies that the functions $\lambda_z, h_z, \nu_z$ exend to holomorphic functions in a neighbourhood 
$J$ of $[0,1)$ so that 
\be\label{complexeigen}
\lambda_z\ne 0, \quad {\M}_zh_z=\l_zh_z, \quad {\M}_z^*\nu_z=\l_z\nu_z, \quad
\nu_z(h_z)=1
\ee
for $z\in J\cup\{1\}$. Given $f\in \Ft$ and $z\in J\cup\{1\}$ 
we decompose
\be \label{decomp1}
{\M}_z\, f = \l_z\, \nu_z (f) \, h_z + {\Nc}_z\, f.
\ee
where the subspace generated by $h_z$ is one-dimensional, whereas
${\Nc}_z$ maps $\Ft$ onto the subspace $\{f\in \Ft\, : \, \nu_z(f)=0\}$. 
Moreover its iterates can be written as ${\Nc}_z^n={\M}_z^n{\cal Q}_z$ where
${\cal Q}_z$ is the spectral projection
valued function $1-{\cal P}_z$ with ${\cal P}_z:=h_z\cdot \nu_z$. 
Now, for each $z\in J\cup\{1\}$ one has the decomposition $\Ft =\Ft^{(1,z)}\oplus \Ft^{(2,z)}$ with 
$\Ft^{(1,z)}={\cal P}_z \Ft$ and $\Ft^{(2,z)}={\cal Q}_z \Ft$. On the other hand, since $h_z$ and 
$\nu_z$ are holomorphic there exists a bounded operator-valued function
${\cal U}_z:\Ft \to \Ft$ with the property that the inverse  
${\cal U}^{-1}_z$ exists as a bounded operator on $\Ft$ and both 
${\cal U}_z$ and ${\cal U}^{-1}_z$
are holomorphic in $J$
and satisfy  (see [Ka], Subsection 7.1.3)
\be\label{anarel}
 h_z = {\cal U}_z\, h,\qquad \nu_z=\nu \,\,{\cal U}_z^{-1}.
\ee
This entails that the pair $\Ft^{(1,1)}$ and $\Ft^{(2,1)}$
decomposes the operator ${\hat{\M}}_{z}:={\cal U}^{-1}_z{\M}_{z}{\cal U}_z$
for all $z$ in $J\cup\{1\}$, and the eigenvalue problem for the part
of ${\M}_{z}$ in $\Ft^{(1,z)}$ is equivalent to that for the part of ${\hat{\M}}_{z}$
in $\Ft^{(1,1)}$ (the eigenvalue being $\l_z$ in both cases).
In particular, it is easily seen that ${\cal U}_z1={\cal U}^{-1}_z1=1$.

\noindent
Now, since ${\Nc}_z$ is holomorphic in $J$, its spectral radius is a lower semicontinuous
function of $z$. Therefore,  since for $z\in (0,1]$ the spectral radius
of ${\Nc}_z$ is strictly smaller than $\l_{z}$,  there is a neighbourhood $H$ of $z=1$ and $\epsilon >0$ such that 
the spectral radius
of ${\Nc}_z$ is smaller than $1-2\epsilon$ for all $z\in H\cap J$. Spectral radius formula 
then implies that $\Vert {\Nc}_z^n \Vert_{\t}\leq (1-\epsilon)^n$ for $n$ large enough and thus
$z\to (1-{\Nc}_z)^{-1}$ is a holomorphic operator-valued function of $z\in H\cap J$.
The above discussion along with Theorem \ref{spectrum} and (\ref{decomp1}) yield the following result.
\begin{proposition} \label{specdet2} 
The function $z\to (1-{\M}_z)^{-1}$ is holomorphic 
in $\D$ and extends continuously to
$\Dc \setminus \{1\}$. In particular, for each $f\in \Ft$ and $z$ in some neighbourhood of $z=1$ (with $z\ne 1$),
we have
$$
(1-{\M}_z)^{-1}f = Q(z)\, h_z \, \nu_z(f) +  (1-{\Nc}_z)^{-1}f
$$
where $Q(z)=(1-\lambda_z)^{-1}$ and $(1-{\Nc}_z)^{-1}f$ is continuous at $z=1$.
\end{proposition}

\section{Transfer operators, invariant measures, ergodic degree and more} \label{section4}
We now go back to the original space $\O$. We shall denote by $\Ft(\O_0)$ the lift of $\Ft(\S)$
with the map $\iota$ (or simply $\Ft$ whenever the underlying space is clear). For notational simplicity'sake, 
we shall use the same symbols $\psi$, $h_z$, $h$, $\nu_z$, $\nu$, $\rho$ to denote 
the lift of the corresponding objects dealt with in the previous Section, as well as
${\M}_z:\Ft(\O_0)\to \Ft(\O_0)$.

\noindent
Now, the transfer operator ${\cal L}: {\cal C}(\O) \to {\cal C}(\O)$
associated to the shift map $T$ on $\O$
and the weight function $\p(\o):=\exp V(\o)$, is defined by
\begin{eqnarray}\label{operatorP}
\left({\cal L} g\right) (\o ) &=& \sum_{T(\o')=\o} \p(\o')\, g(\o')\nonumber \\
&=&\p(0\o)\, g(0\o )+\p(1\o)\, g(1\, \o )=:({\L}_0 + {\L}_1 )\, g(\o ). 
\end{eqnarray}
From Property 4 it follows that
$r({\cal L}) = 1$. 
Moreover we can write,
\begin{eqnarray}\label{chain}
\left({\M}_{z}f\right)(\o) &=& 
\sum_{k=1}^{\infty} z^{k}\,\psi
(0^{k-1}1\,\,\o)\cdot f(0^{k-1}1\,\,\o)\nonumber \\
&=&\sum_{k=1}^{\infty} z^{k}\, \prod_{i=0}^{k-1} \p
(0^{k-i-1}1\,\,\o)\cdot f(0^{k-1}1\,\,\o)\nonumber \\
&=&\sum_{k=1}^{\infty} z^{k}\,{\cal L}^{k}(f\cdot 1_{A_k})(\o) \\
&=&\sum_{k=1}^{\infty} z^{k}\, \left({\cal L}_1 {\cal L}_0^{k-1}f\right) (\o)\nonumber \\
&=&z\,{\cal L}_1 (1-z {\cal L}_0)^{-1}f(\o)\nonumber 
\end{eqnarray}
where (\ref{inducedversion}) and (\ref{operatorP}) have been used in the second and third equalities.
From (\ref{chain}) we obtain the following algebraic relation 
between ${\M}_{z}$ and
${\cal L}$ (see also \cite{HI} where a similar relation has been exploited to study statistical properties of rational maps):
\begin{proposition} \label{identity} For $z\in \Dc$
and for any $f$ such that $(1-{\L}_0)f\in \Ft$ we have
$$
(\, 1-{\M}_{z}\,  )\, (\, 1-z\, {\L}_0 \, )\, f = (1-z\, {\L})\, f.
$$
\end{proposition}
From this identity and Proposition \ref{specdet2} we get
\begin{corollary}\label{decco} The function $z\to (1-z{\L})^{-1}$ is holomorphic 
in $\D$ and extends continuously to
$\Dc \setminus \{1\}$. For all
$f$ s.t. $(1-{\L}_0)f\in \Ft$ and $z$ in some neighbourhood of $z=1$ (with $z\ne 1$) we have,
\be
(1-z\, {\L})^{-1}f = Q(z)\, e_z \, \nu_z(f) +  (1-z{\L}_0)^{-1}(1-{\Nc}_z)^{-1}f,
\ee
where $e_z=(1-z\,{\L}_0)^{-1}h_z$.
\end{corollary}
Other consequences of Proposition \ref{identity} are the following:
\begin{proposition} \label{relations} Let $z\in \Dc \setminus \{0\}$. Then 
$1$ is an eigenvalue of ${\M}_{z}$
if and only if $1/z$ is an eigenvalue of
${\L}$ and they have the same geometric 
multiplicity. 
Furthermore, the
corresponding eigenfunctions $f_z$ of ${\L}$ and $g_z$ of 
${\M}_{z}$ are related by $g_z = (1-z{\L}_0) f_z$ or else  
$f_z = \sum_{k=0}^{\infty} z^k{\L}_0^k g_z$.
\end{proposition}
{\it Proof.} Assume that ${\M}_{z}g_z = g_z$. From
Proposition \ref{identity} it follows that 
$(1-z{\L})\sum_{k=0}^{\infty} z^k{\L}_0^k g_z = 0$. Conversely, 
assume that $z{\L} f_z = f_z$, then 
$(1-{\M}_{z})(1-z{\L}_0)f_z=0$. $\qed$
\vskip 0.1cm 
\noindent
\begin{proposition} \label{measures} The $\s$-finite measure $\mu = e\cdot \nu$ with 
$$
e = \sum_{k=0}^{\infty} {\L}_0^k h\quad\hbox{or else}\quad h=(1-{\L}_0)e={\L}_1e
$$
(here ${\M}h=h$ and thus ${\L}e=e$) is 
$T$-invariant. To any Borel
subset $E$ of $\O_0$ it assigns the weight
$$
\mu (E) = \sum_{k\geq 0} \rho \left(T^{-k}E\cap D_k\right)
$$
where $D_k=\{\tau > k\}= \cup_{l>k}A_l $. 
In particular $\mu (A_k) =\rho (D_{k-1})$. 
\end{proposition} 
Conversely, the measure $\rho$ is obtained by pushing backward $\mu$ with the map $T_1:\O \to \O$
given by $T_1\o = 1\o$, i.e.
\be\label{converse}
\rho (E) = (\mu\circ T_1)(E).
\ee
In particular we have
\be
\rho(A_n)=\mu(B_n), 
\ee
with
\be
B_n=T_1(A_n)=\{\o\in A_1 \,:\, \min\{k\geq 1 \, :Ê\, T^k\o \in A_1\}=n\}.
\ee
Moreover, we have the following chain of formal identities:
\be \label{Kac}
\mu (\O) = \sum_k \mu (A_k) =\sum_k \rho (D_{k-1}) = 
\sum_k  \sum_{l\geq k} \rho (A_l) =\sum_k k\cdot  \rho (A_k) =\rho (\tau)
=:M_1
\ee
Here $M_1$ denotes the mean first passage time in the state $1$ for the dynamical system
$(\O,T,\mu)$ (which might be infinite).
Thus (\ref{Kac}) is a version of Kac's formula. More generally, let $M_{\gamma}$, $\gamma \geq 0$, 
be the family of moments defined by
\be
M_{\gamma}:=\sum_k k^\gamma\cdot  \rho (A_k).
\ee
\begin{definition} \label{degree} The triple $(\O,T,\mu)$ is said to have {\rm ergodic degree} $d$
if $M_{d+1}=\infty$ but $M_{d+1-\epsilon}<\infty$, $\forall \epsilon >0$.
\end{definition}
\begin{remark}
{\rm  One could also define the ergodic degree using moments wrt the $\s$-finite measure $\mu$
by saying that $(\O,T,\mu)$ has {\rm ergodic degree} $d$
if $\sum_k k^d\cdot  \mu (A_k)=\infty$ but 
$\sum_k k^{d-\epsilon}\cdot  \mu (A_k)<\infty$, $\forall \epsilon >0$.
The above definition
may appear somewhat strict. In particular it may happen that setting
${\tilde d}=\inf\{\gamma : M_{\gamma+1}=\infty\}$ one has $\sum_k k^{\tilde d}\cdot  \mu (A_k)<\infty$. 
In this case the ergodic degree does not esist. Take for instance
$\mu(A_k)=k^{-1}(\log k)^{-2}$. Then we have $M_{\gamma+1} =\infty$
for all $\gamma >0$ (so that  ${\tilde d}=0$) but $\mu(\O)=M_1<\infty$.
\noindent
On the other hand, the above definition has the advantage to make a neat distinction between
finite measure case $\mu (\O)<\infty$, where $d>0$, and the infinite one, where $d\leq 0$. 
In particular, notice that $M_{0} =1$, 
so that the ergodic degree satisfies $d>-1$, provided it exists.
For the Markov chain of Example 1, this notion is related to
some already used in the literature (see, e.g., [Is2] and references therein).
In particular, if $-1<d \leq 0$ one has a null-recurrent chain, whereas for
$d> 0$ it is positive recurrent. The Markov measure $\mu$ is given by 
$\mu([x_0\cdots x_n])=\pi_{x_0} p_{x_0x_1}\cdots p_{x_{n-1}x_n}$ with 
$\pi_{k}=\sum_{\ell \geq k}p_\ell$ and is infinite if $d\leq 0$.}
\end{remark}

\begin{remark} {\rm 
By
Property 3 we have
\be
e^{-C_3\theta}  \leq {\psi(0^{n-1}1\o) \over\psi(0^{n-1}1\o')}\leq e^{C_3\theta}.
\ee
Moreover, using normalization and conformality of the measure $\nu$ we have 
\begin{eqnarray}
\nu (A_n) &=&\int_\O 1_{A_n}(\o) d\nu (\o) =\int_\O {\M} 1_{A_n}(\o) d\nu (\o)\nonumber \\ 
&=&\int_\O  \psi(0^{n-1}1\o)  d\nu (\o) = \psi(0^{n-1}1\o^*) 
\end{eqnarray}
for some $\o^* \in \O$. Since $h\asymp 1$ we
therefore have\footnote{{\sl Notational warning}: Here and in the sequel, 
for two sequences $a_n$ and $b_n$ we shall write:
\begin{itemize} 
\item $a_n \approx b_n$ if the ratio $a_n/b_n$
grows slower than any power of $n$, or decays slower than any inverse power of $n$, as $n\to \infty$; 
\item $a_n \asymp b_n$, if
$C^{-1} \leq a_n/b_n \leq C$ for all $n$ and fixed $C\geq 1$;
\item $a_n \sim b_n$
if the quotient $a_n/b_n$ tends to unity as $n\to \infty$. 
\end{itemize}}
\be\label{asymp}
\rho (A_n) \asymp \psi(0^{n-1}1\o)
\ee
uniformly in $\o \in \O$. 
Suppose now that $\psi(0^{n-1}1\o)\asymp n^{-\alpha}$,
for some $\alpha >1$. Then one finds $d=\alpha -2$. The case 
$\psi(0^{n-1}1\o)\asymp n^{-\alpha}\,L(n)$
where $L(n)$ a function slowly varying at infinity, 
(i.e. $L(cn) \sim L(n)$ for every positive $c$) is more delicate. If for instance $L(n)=\log n$
then again $d=\alpha-2$. If instead $L(n)=(\log n)^2$ then, as already noted, $d$ does not exist because
$\sum n^{d+1-\alpha}(\log n)^{-2}<\infty$ for $d=\alpha-2$ although $\alpha-2=\inf\{d : M_{d+1}=\infty\}$.}
\end{remark}
\begin{remark} {\rm For every finite $d$ the fixed function $e$ 
for the operator ${\L}$ 
extends to a unique extended-real-valued
function on $\O$, still denoted by $e$, 
such that $e$ is finite except at $0^{\infty}$
where it takes the value $+\infty$. Moreover it is easy to check that for each $k$ there is constant
$D=D(k)$ so that
$\Var_n \, (e\cdot 1_{A_k}) \leq D\, \t^{n}$. On the other hand the function $h=(1-{\L}_0)e$
is uniformly bounded on $\O$ and $\Var_n h \leq C\, \t^{n}$. 
Finally, we may extend the measures $\rho$ and $\mu$ to the whole space $\O$ by putting 
$\rho(\{0^\infty\})=\mu(\{0^\infty\})=0$. }
\end{remark}
\begin{remark}\label{weakgibbs} {\rm As shown in the previous Section (see Corollary \ref{gibbs}), the probability measure $\rho$
can be viewed a Gibbs measure on $\O_0$ which is invariant under the action of
the induced shift $T^\tau$. Now, if $d>0$ 
the corresponding property for the $T$-invariant probability measure ${\hat \mu}:=\mu/M_1$
would be
\be
{{\hat \mu} ([\o_0\cdots \o_{n-1}])\over \exp{\sum_{k=0}^{n-1}V(T^k\o)}}\asymp 1.
\ee
On the other hand, if $d$ is positive but finite, 
this is not the case. Indeed, 
consider the cylinder set $A_n$ (see (\ref{cylinder})).
According to Proposition \ref{measures} and (\ref{asymp}) we have
\be
\mu ([0^{n-1}1]) \asymp \sum_{l\geq n}\psi(0^{l-1}1\o)
\ee
uniformly in $\o \in \O$.
Now, if $(\O,T,\mu)$ has finite ergodic degree then for all $\o\in \O$,
\be 
\psi (0^{n-1}1\o) =o\left(\sum_{l\geq n}\psi(0^{l-1}1\o) \right),
\ee
and therefore
\be
{{\hat \mu} ([0^{n-1}1])\over \exp{\sum_{k=0}^{n-1}V(0^{n-k-1}1\o})}\asymp
\Delta (n)
\ee
with $\Delta (n) \approx n$. This can be interpreted as a {\it weak Gibbs property}
\cite{MRTVV}, \cite{Yu}. It appears to be related to the fact that the potential $V$ has (exactly) two
equilibrium states (in the sense of Walters). Let us see this point in more detail. According to \cite{Wal2}, and using Property 4 of the potential function $V$, a
$T$-invariant probability measure $m$ on $\O$ is called an {\sl equilibrium state} for $V$ if it
satisfies
\be\label{variational}
h_m(T)+m(V)=0,
\ee
where $h_m(T)$ denotes the measure-theoretic entropy of $T$ wrt $m$. Now, on the one hand, by Property 1 it is plain that
the point measure $\delta_{0^\infty}$ concentrated at $\{0^\infty\}$ satisfies $h_{\delta_{0^\infty}}(T)+
\delta_{0^\infty}(V)=0$.
On the other hand, if $M_1<\infty$,
it follows from the results of Section \ref{opvalpowser} and \cite{Wal2}
that the measure $\rho$ is the {\sl unique} equilibrium state
for $W$, namely $h_\rho(T^\tau)+\rho(W)=0$ and any other $T$-invariant (and thus also $T^\tau$-invariant) 
measure ${\tilde \rho}$ 
(with ${\tilde \rho}(\{0^\infty\})=0$)
satisfies
$h_{\tilde \rho}(T^\tau)+{\tilde \rho}(W)<0$. Moreover we have
\begin{eqnarray}
\mu(V) &=& \int_{\O} V(\o)\, e(\o) \nu(d\o) =\int_{\O} V(\o)\, \sum_{n=0}^{\infty}{\L}_0^n h(\o)\; \nu(d\o)\nonumber \\
&=& \sum_{n=0}^\infty \int_{D_n} V(T^n\o)\, h(\o)\, \nu(d\o)=
\sum_{n=1}^\infty \int_{A_n} \left(\sum_{k=0}^{n-1}V(T^k\o)\right)\, h(\o)\, \nu(d\o)\nonumber\\
&=&\int_{\O} W(\o)\, h(\o) \, \nu(d\o) = \rho (W).
\end{eqnarray}
In addition, taking natural extension and using Abramov's formula \cite{Ab} we get 
$h_\rho(T^\tau)=M_1 \cdot h_{\hat \mu}(T)$ so that, finally,
\be
h_{\hat \mu}(T)+ {\hat \mu}(V) ={h_\rho(T^\tau)+\rho(W)\over M_1}=0.
\ee
Whence, any $T$-invariant measure $m$ which satisfies (\ref{variational}) is a convex combination of
$\delta_{0^\infty}$ and ${\hat \mu}$ (see also \cite{Ho} and \cite{FL} for related examples). 
}
\end{remark}

\section{The pressure function $P(z)$} \label{pressione}
We now characterize the 
behaviour 
of the pressure function $P(z)$ and of its derivatives when $z\uparrow 1$ (and thus of the leading eigenvalue $\lambda_z =e^{P(z)}$) 
in terms of suitable expectations wrt the equilibrium state $\rho$. 
We first introduce some further quantities. For $\gamma \geq 0$ and $|z|<1$, set
\be
\label{moments}
M_\gamma (z) := \sum_{k=1}^\infty z^{k}\, k^\gamma \, \rho (A_k)
\ee
so that $M_\gamma (1)\equiv M_\gamma$ (provided it exists).
Notice that $M_0(z)=\nu ({\M}_z h)$. More generally define, for $m\geq 0$,
\be \label{derivatives}
{\M}_z^{(m)}f (\o) :=\sum_{k=1}^{\infty} z^{k}\,k^m\, 
\left({\cal L}_1 {\cal L}_0^{k-1}f\right) (\o) = {\M}_z(f\cdot \tau^m) (\o) ,
\ee
so that ${\M}_z^{(0)}\equiv {\M}_z$.
Reasoning as in the proof of Lemma \ref{radius} one sees that if $(\O,T,\mu)$ has ergodic degree $d>0$
then, for all integers $m$ with $0\leq m< d+1$, the power series of 
${\M}_{z}^{(m)}$ when acting on ${\Ft}$ has radius of convergence
bounded from below by $1$ and, moreover, 
it converges absolutely at every point of the unit circle. In addition, from (\ref{moments}) and 
(\ref{derivatives})  we have, for $m> 0$,
\be\label{momenti}
M_m (z) =(zD)^m M_0(z)= \nu\left( {\M}_z^{(m)} h \right),
\ee
where $zD$ denotes the differential operator
$z(d/dz)$. It is not difficult to realize that, whenever it is defined for $0\leq z\leq 1$, 
the operator ${\M}_{z}^{(m)}$ has leading eigenvalue $\lambda_z^{(m)}=(zD)^m\lambda_z$.

\begin{theorem} \label{analytic} The function $z\to \lambda_z\equiv e^{P(z)}$ is 
analytic in a complex open neighbourhood of $[0,1)$. 
Moreover, we have the following properties:
\begin{enumerate}
\item If $-1<d\leq 0$ then $(zD)\lambda_z \sim M_1(z)$ as $z\uparrow 1$;
\item if $d>0$ then $P(z)\in C^1((0,1])$ and 
$(zD)P(z)\Big|_{z=1}=M_1=\rho(\s_0)$;
\item if $d>1$, then $P(z)\in C^2((0,1])$ and 
$$
(zD)^2P(z)\Big|_{z=1}=\sigma^2:=\sum_{n\geq 0} \bigl[ 
\rho(\s_0\s_n)-\rho(\s_0)\rho(\s_0)\bigr];
$$
\item more generally, if
$m-1< d \leq m$ for some $m\in \Z_+$ then $P(z)\in C^m((0,1])$ and,  for $1\leq \ell \leq m$,
$$
(zD)^\ell P(z)\Big|_{z=1}=\sum_{n_1\geq 0}\cdots \sum_{n_{\ell -1}\geq 0}
U_\ell(\s_0,\s_{n_1},\dots ,\s_{n_{\ell-1}})\, < \, \infty,
$$
where 
$$
{ U}_{\ell}(\s_{k_1},\dots ,\s_{k_\ell}):=
{\partial^\ell \over \partial t_1 \dots \partial t_\ell}
\log \rho
\left(\exp (\sum_{i=1}^\ell t_{i}\s_{k_i})\right)\Big|_{t_1=\cdots =t_\ell=0}.
$$
In addition,
$(zD)^{m+1}\lambda_z - M_{m+1}(z)={\cal O}(1)$ as $z\uparrow 1$.
\end{enumerate}
\end{theorem}
{\it Proof.} In this proof we shall use the notation of Section \ref{opvalpowser} and consider ${\M}_z$
as acting on $\Ft(\S)$. We first notice that, 
as a consequence of Lemma \ref{radius}, $z\to {\M}_{z}$ is an analytic family in the sense of Kato
for $z$ in the open unit disk. 
Theorem \ref{spectrum} and Kato-Rellich Theorem (see, e.g., Thm. XII.8 in \cite{RS}) then imply that
$\exp P(z)$ extends analytically in a complex open neighbourhood of $[0,1)$.
In addition, since $\Lambda_n(z)$ is monotonically increasing for $0\leq z\leq 1$, 
so is $\exp P(z)$, and therefore $\exp P(z) < \exp P(1)=1$ for $z\in [0,1)$. This proves
the first statement.
Now, from the proof of Theorem \ref{spectrum} we know that for $z$ real and positive
$\lambda_{z,N}=\exp P_N(z)$ is the (simple) eigenvalue with largest modulus
of ${\M}_{z,N}:{\Ft}(\S_N)\to {\Ft}(\S_N)$, 
where $P_N(z)$ is defined in (\ref{P_N}). We now define a function $K_N(z)$ by setting 
\be
\label{decomposition_N}
\exp P_{N}(z)=\exp P_N(1)\left(\sum_{k=1}^Nz^k\, \rho_N(A_k) + K_{N}(z) \right)
\ee
where $\rho_N=h_{1,N}\cdot \nu_{1,N}$ is
the equilibrium state on $\S_N$ for the function $W$ and 
$\rho_N(A_k) = \rho_N \{\s \in \S_N : \s_0=k\}$. Clearly we have $K_{N}(0)=K_{N}(1)=0$.
Furthermore, 
from the proof of Theorem \ref{spectrum}
it follows that for all $0\leq z \leq 1$
the function $\exp P_N(z)$ converges to $\exp P(z)$ as $N\to \infty$.
We also know from (\cite{CI}, Theorem 2.1) that
$\exp P_N(1)\to 1$ and $\rho_N \to \rho$
uniformly as $N\to \infty$, where $\rho=h\cdot \nu$ 
is the equilibrium state on $\S$
for the function $W$. Therefore, for $0\leq z \leq 1$ and $N\to \infty$ we have that
$K_{N}(z)$ tends pointwise to a function $K(z)$ so that 
\be
\label{decomposition}
\exp P(z) =  \sum_{k=1}^{\infty}z^k \rho (A_k) +K(z),
\ee
The function $K(z)$ is analytic in a complex open neighbourhood
of $[0,1)$ and satisfies $K(0)=K(1)=0$.
We now
substitute $z=e^t$ 
into (\ref{decomposition_N}) and expand the sum $\sum_{k=1}^Nz^k\, \rho_N(A_k)$ in powers of $t$:
\be \label{powersoft}
\exp P_N(e^t)=\exp P_N(1)\left(
\sum_{m=0}^{\infty} {\rho_N(\s_0^m)\over m!}t^m + K_{N}(e^t)   \right) 
\ee
where $\rho_N(\s_0^m) =\sum_{k=1}^{N} k^m\, \rho_N(A_k)$.
Clearly $(zD)^\ell ={d^\ell / dt^\ell}$ for $z=e^t$.
Now observe that
$P_N(e^t)$ is the pressure of 
the function $W+t\, \s_0$ restricted to $\S_N$.
It is a standard result in the theory of equilibrium states (see [Ru4], Chapter 5)
that, under the conditions assumed here, $P_N(e^t)$ is analytic 
in some neighbourhood of $t=0$ and its derivatives at $t=0$ are
related to suitable moments of $\rho_N$ 
(see \cite{Ru4}, p.99-100). 
In particular, we have
\be \label{id1}
{dP_N(e^t)\over dt}\Big|_{t=0}=\rho_{N}(\s_0),
\ee
and
\be \label{id2}
{d^2P_N(e^t)\over dt^2}\Big|_{t=0}=\sum_{n\geq 0} \bigl[ \rho_{N}(\s_0\s_n)-
\rho_{N}(\s_0)\rho_{N}(\s_0) \bigr].
\ee
Now, differentiating twice (\ref{powersoft}) at $t=0$, and using (\ref{id1}) and (\ref{id2}), we get
$$
{dK_N(e^t)\over dt}\Big|_{t=0}=0\quad\hbox{and}\quad
{d^2K_N(e^t)\over dt^2}\Big|_{t=0}= \sum_{n> 0} \bigl[ \rho_N(\s_0\s_n)-
\rho_N(\s_0)\rho_N(\s_0) \bigr].
$$
Notice that, for any fixed $N$, the expression in square brackets
decreases exponentially with $k$. More specifically, set ${\M}_{N}\equiv {\M}_{1,N}$, 
$\lambda_N\equiv \lambda_{1,N}$, ${h_N}\equiv h_{1,N}$, 
${\nu_N}\equiv \nu_{1,N}$ where $\lambda_{1,N}, h_{1,N}, \nu_{1,N}$
are as above, and define
\be
X_0^N= {{\M}_{N}\left( \, {h_N}\, \s_0\, \right)\over 
\lambda_N\, \rho_N(\s_0) }.
\ee 
Reasoning as in the proof
of Lemma \ref{radius}, one can find
a positive constant $C$, independent of $N$, such that,
for $N$ large enough, we have
$\Vert X_0^N \Vert_{\t} < C$.
On the other hand, from the proof of Theorem \ref{spectrum} it follows 
that there exist two constants $M>0$ and $0<\vartheta <1$, independent of $N$, 
such that, for $N$ large enough and for all $v\in {\Ft}(\S_N)$,
we have 
$$
\Vert \lambda_N^{-k}{\M}^{k}_{N}v -{h_N}\cdot \nu_N(v) \Vert_{\t}
\leq \, M\, \vartheta^k \, \Vert v \Vert_{\t}.
$$
Whence, taking $v=X_0^N$ and using ${\nu_N}(X_0^N)=1$, we get
$$
\Vert \lambda_N^{-k}{\M}^{k}_{N}X_0^N -{h_N} \Vert_{\t}
\leq \, C \, M\, \vartheta^k
$$
and therefore, for $k>0$, 
$$
|\rho_N(\s_0\s_k)-\rho_N(\s_0)\rho_N(\s_0)| 
\leq C \, M\, \vartheta^{k-1} \, \rho_N(\s_0)\rho_N(\s_0).
$$
This gives 
\be
{d^2K_{N}(e^t)\over dt^2}\Big|_{t=0}\leq 
{C \, M\, \over 1-\vartheta} \, \rho_N(\s_0)\rho_N(\s_0).
\ee
From this discussion we obtain that the 
following properties hold uniformly in $N$:
for all $d>-1$, we have
\be
\lim_{t\to 0_-}{dK_{N}(e^t)\over dt} =0
\ee
and moreover, if $d>0$, then 
\be
{dK_{N}(e^t)\over dt}\Big|_{t=0} =0\quad\hbox{and}\quad
{d^2K_{N}(e^t)\over dt^2}\Big|_{t=0} < \infty.
\ee
We can actually say more.
First observe that if $d\leq 0$ then $\rho_N(\s_0)\to \infty$ as 
$N\to \infty$. 
In addition, if $m-1< d \leq m$, with $m>0$, then
the expectations
$\rho_N(\s_0^{\ell})$ are
bounded uniformly in $N$ for $1\leq \ell \leq m$, 
but $\rho_N(\s_0^{m+1})\to \infty$ as $N\to \infty$. 
On the other hand, if $d>m-1$ the derivatives 
${d^{{\ell}}K_{N}(e^t)/ dt^{{\ell}}}|_{t=0}$ are
bounded uniformly in $N$ for $1\leq {\ell} \leq m+1$, as we are now going to show
(notice however that the case $m=1$ has already been discussed).

\noindent
Using induction, one can compute
the $\ell$-th derivative (with $\ell >1$) of the pressure $P_N(e^t)$ at $t=0$ as:
\be \label{id3}
{d^\ell P_N(e^t)\over dt^\ell}\Big|_{t=0}=\sum_{k_1\geq 0}\dots \sum_{k_{\ell-1}\geq 0}
{ U}_{\ell,N}(\s_0,\s_{k_1},\dots ,\s_{k_{\ell-1}}),
\ee
where ${ U}_{\ell,N}$ is 
given by
$$
{ U}_{\ell,N}(\s_{k_1},\dots ,\s_{k_\ell})=
{\partial^\ell \over \partial t_1 \dots \partial t_\ell}
\log \rho_N\left(\exp (\sum_{i=1}^\ell t_{i}\s_{k_i})\right)\Big|_{t_1=\cdots =t_\ell=0}.
$$
The function ${ U}_{\ell,N}$ can be considered as a particular 
version of what in statistical mechanics is called the $\ell$-th Ursell function 
(see, e.g., \cite{Si}, Section II.12).
Let us write
\be \label{id4}
{d^{m+1}P_N(e^t)\over dt^{m+1}}\Big|_{t=0}= { U}_{m+1,N}(\s_0,\dots ,\s_0) + { V}_{m+1,N}.
\ee
From the characterization of the
${ U}_{\ell,N}$'s given in (\cite{Si}, Corollary II.12.7) we have
that, for any $\ell < m+1 $, ${ U}_{\ell,N}(\s_0,\s_{k_1},\dots ,\s_{k_{\ell-1}})$ 
can be written as a 
linear combination of products 
$$
\rho_N(\s_0^{r_1})\dots \rho_N(\s_0^{r_k})\,
\bigl( \, \rho_N(X_1X_2)-\rho_N(X_1)\rho_N(X_2) \, \bigr)
$$
for suitable $r_1,\dots,r_k$ and with $X_1$ a product
of functions from among $\s_0,\dots,\s_{k_{\ell-1}}$ and $X_2$ 
a product of functions from among $\s_{k_{\ell}},\dots,\s_{k_{m}}$.
Thus, reasoning 
as above, one obtains that, for $N$ large enough, 
$|{ V}_{m+1,N}|$ is bounded from above by a linear combination, whose
coefficients depend of $m+1$ but not of $N$, of products
$\rho_N(\s_0^{r_1})\dots \rho_N(\s_0^{r_n})$,
with $r_1+\cdots + r_n=m+1$ and $m+1\neq r_i\geq 0$ for all $i\in \{1,\dots ,n\}$.
Notice that the leading term (with respect to the limit $N\to \infty$) is
given by
$\rho_N(\s^{m}_0)\,\rho_N(\s_0)$ and corresponds to the choice
$k_1=k_2=\dots =k_{m}\neq 0$ in (\ref{id3}).
Moreover, differentiating (\ref{powersoft}) $m$ times at $t=0$ and using (\ref{id3}),
along with the combinatorial features of the $U_{\ell,N}$'s
(see \cite{Si}, Section II.12),
one sees, inductively in $\ell$, that
${d^{m+1} K_{N}(e^t)/ dt^{m+1}}$ at $t=0$ is equal to $V_{m+1,N}$ 
plus a linear combination of the
${ V}_{k,N}$'s with $k=2,\dots ,m$. This finishes the proof. $\qed$

\vsni
\noindent
\section{A Markov approximation}\label{markovapprox}
We now consider the {\sl renewal sequence} $a_0,a_1,\dots $ 
associated with the sequence
$p_1,p_2,\dots$, with $p_n \equiv\rho (A_n)$, which is generated by the recurrence relation:
\be\label{renew}
a_0=1\quad\hbox{and}\quad a_n=p_n+a_1\, p_{n-1}\cdots + a_{n-1}\, p_1
\quad\hbox{for}\quad n\geq 1.
\ee 
Its generating function $A(z)$ is given by
\be \label{generatingfct1}
A(z)=\sum_{n=0}^{\infty} a_n z^n= \left(1-\sum_{n=1}^{\infty}p_nz^n\right)^{-1}=
\left( (1-z)\sum_{n=0}^{\infty}d_nz^n\right)^{-1},
\ee
where we have written $d_n =\sum_{k >n}p_k= \rho(D_n)$ and the $D_n$'s are defined in Proposition \ref{measures}.
If $d>0$, another generating function of interest is that of the numbers $b_n:=M_1\, a_n-1$, that is
\be \label{generatingfct2}
B(z)=\sum_{n=0}^{\infty} b_n z^n= M_1\, A(z)-{1\over 1-z}={D^{(1)}(z)\over D(z)},
\ee
where we have put
\be \label{D^1}
D^{(1)}(z) = \sum_{n\geq 0} d^{(1)}_n z^n,\qquad d^{(1)}_n := \sum_{l>n}d_l = \mu (D_{n+1})
=\mu(\tau>n+1),
\ee
and in the last identity we have used Proposition \ref{measures}. The asymptotic behaviour of the sequences $a_n$ and $b_n$ is described in
the following result, whose proof can be readily extracted from
(\cite{Is2}, Theorem 2 and Remark 2).
\begin{proposition} \label{generfcts} Suppose the ergodic degree $d$ exists and is finite. Then the
power series expansions of (\ref{generatingfct1}) and (\ref{generatingfct2}) define holomorphic functions $A(z)$ and $B(z)$
in $\D$ which converge uniformly at every point of the unit
circle with the exception of $z=1$, where they have a non-polar singular point.
Moreover, one has
the following asymptotic behaviour 
of their coefficients: for every $d> -1$ we have 
$$
a_n \to {1\over M_1}\quad\hbox{as}\quad n\to \infty, 
$$
where $1/M_1\equiv 0$ if $M_1 = \infty$. Furthermore, if $-1<d\leq 0$ we have
$$
a_n \approx n^d,
$$
whereas, if $d>0$, 
$$
b_n\equiv M_1\, a_n-1\sim \left[ M_1^{-1} \mu (\tau > n) \right] \approx n^{-d}.
$$
\end{proposition}

\noindent
Now, if we view the partition set $A_n$ as the $n$-th `state' 
for the dynamical system $(\O,T,\mu)$, the quantity $p_n\equiv \rho (A_n)$ can be interpreted 
as the $\rho$-probability 
that a first passage in the state $1$ occurs after $n$ iterates.
We may then consider the quantity 
\be\label{u_n}
u_n:=\mu (A_1\cap T^{-n}A_1),
\ee 
that is the $\mu$-probability 
to observe a {\sl return} in the state $1$ after 
$n$ iterates (recall that $\mu (A_1)=1$). 
Using (\ref{converse}) backward we have 
\be
u_n=\mu ( T_1 (T^{-(n-1)} A_1))=\rho \left( T^{-(n-1)} A_1\right)
\equiv \rho \left( T^{n-1}(x)\in A_1\right),
\ee
the last quantity being the
$\rho$-probability to observe
a {\sl passage} in the state $1$ after $n-1$ iterates (for the first time or not).
Another interpretation of $u_n$ is the following.
Let $N_n(\o):=1_{A_1}(\o)+\cdots +1_{A_1}(T^{n-1}(\o))$
be the number of passages in the state $1$ up to the $n$-th iterate of the map $T$.
Let also $s_n(\o)=\s_0(\o)+\s_1(\o)+\cdots +\s_{n-1}(\o)$ 
be the total number of iterates of $T$ needed to observe $n$ passages 
in the state $1$. Then notice that 
$(N_n=k)=(s_k\leq n< s_{k+1})=(s_k\leq n)-(s_{k+1}\leq n)$. Thus
$\rho (N_n=k)=\rho (s_k\leq n)-\rho (s_{k+1}\leq n)$, which is the same as 
$\rho (s_k\leq n)=\sum_{r=k}^n\rho (N_n=r)$. Moreover
$\rho (s_k=n)=\rho (s_k\leq n)-\rho (s_k\leq n-1)$ for $k<n$ and $\rho (s_n=n)=\rho (s_n\leq n)$. 
Therefore we have the following chain of identities
\begin{eqnarray}\label{passages}
u_n&=&\sum_{k=1}^n \rho (s_k=n) =\sum_{k=1}^n \rho (s_k\leq n)-\sum_{k=1}^{n-1}\rho (s_k\leq n-1)
\nonumber\\
&=&\sum_{k=1}^n \sum_{r=k}^n \rho (N_n=r)-\sum_{k=1}^{n-1}\sum_{r=k}^{n-1}\rho (N_{n-1}=r)
= \rho (N_n)-\rho (N_{n-1}) 
\end{eqnarray}
where $\rho (N_n)$ denotes the mean of the random variable $N_n$ (set $N_0=0$).
Thus, $u_n$ may be regarded as the expected number of passages in the state $1$ 
per iteration of the map $T$ (after $n-1$ iterations), namely as a `renewal density' 
for $(\O,T,\mu)$.

\noindent
A further decomposition using a first passage argument is the following:
\be\label{decomp}
u_n=\sum_{r=1}^{n}\, p_r\, \,
\rho \left(T^{n-1}\o\in A_1\, |\, \o\in A_r \,  \right).
\ee
Now suppose that the sequence of random variables 
$\s_0(\o), \s_1 (\o), \s_2(\o), \dots$ were mutually independent. 
In this case the iteration process
$\{T^n\o\}_{n\geq 0}$ would `start afresh' at each passage in the state $1$ and 
\be\label{indep}
\rho \left(T^{n-1}\o\in A_1\, |\, \o\in A_r  \,  \right)
=\,\rho \left(T^{n-r-1}\o\in A_1\right)=u_{n-r}.
\ee
Direct comparision of (\ref{renew}), (\ref{decomp}) and (\ref{indep}) shows that under the above supposition
we would have $u_n\equiv a_n$. 

\noindent
Furthermore, it has been observed (see, e.g., \cite{Ki}) that any renewal sequence, that is any sequence
generated as in (\ref{renew}) with $p_1,p_2\dots$ satisfying $p_n\geq 0$ and
$\sum p_n \leq 1$, can arise as the diagonal transition probabilites corresponding to a
given state in some Markov chain. In our case, it is easy to check that a Markov chain which 
does the job is that 
discussed in Example 1 with the identification $p_n=\rho(A_n)$.
According to the above discussion,
we shall say that
the Markov chain of Example 1 is the {\it Markov approximation of $(\O,T,\mu)$ with respect to the reference
set $A_1$}. 
Its statistical properties are partially described by Proposition \ref{generfcts}
(see also \cite{Is1}, \cite{Is2}; and \cite{Che} for a related approximation scheme).

\vsni
\noindent
From Theorem \ref{analytic} and Proposition \ref{generfcts} we obtain a first
result about the proximity of $(\O,T,\mu)$ and its Markov approximation.
\begin{corollary} \label{singular} For all finite $d> -1$,
the function $Q(z):=(1-\exp{P(z)})^{-1}$ has a non-polar
singular point at $z=1$ with
$Q(z) \sim A(z)$ as $z\uparrow 1$,
where $A(z)$ is the generating function defined in (\ref{generatingfct1}). 
Moreover we have $Q(z)-A(z)={\cal O}(1)$ and
$(zD)\left[Q(z)-A(z)\right]={\cal O}\left(M_{2}(z)/M_1(z)\right)={\cal O}\left((1-z)^{-1}\right)$ as $z\uparrow 1$.

\noindent
For $d>0$, 
the function $U(z):=M_1\, Q(z) - (1-z)^{-1}$ has the following properties:
\begin{enumerate}
\item if $0<d\leq 1$ then
$U(z)\sim B(z)$ as $z\uparrow 1$,
where $B(z)$ is the generating function defined in (\ref{generatingfct2});
\item if $1<d\leq 2$ we have
$U(z)\to {\di \sigma^2 +M_1^2-M_1\over \di 2M_1}$ and 
$(zD)U(z)\sim (zD)B(z)$ as $z\uparrow 1$;
\item more generally, 
if $m-1<d\leq m$, for $m>0$, then
$(zD)^\ell U(z)$ is uniformly bounded 
for $0\leq \ell <m-1$ but 
$(zD)^{m-1}U(z)\sim (zD)^{m-1}B(z)$ as $z\uparrow 1$. In addition we have 
$(zD)^{m-1}\left[U(z)-B(z)\right]=o\, (1)$ and
$(zD)^{m}\left[U(z)-B(z)\right]={\cal O}(M_{m+1}(z))=o\, \left((1-z)^{-1}\right)$ as $z\uparrow 1$.
\end{enumerate} 
\end{corollary}
{\it Proof.}
From (\ref{decomposition}) we can write
\be\label{dec1}
Q(z)=A(z)\cdot (1+J(z))\quad\hbox{with}\quad
J(z):=\left({A(z)\cdot K(z)\over 1-A(z)\cdot K(z)}\right).
\ee
Now, from the proof of Theorem \ref{analytic} we have that when $z\uparrow 1$ and for all $d>-1$
\begin{eqnarray}\label{kprimo}
(zD)K(z)&=&(zD)\lambda_z - M_1(z)=\lambda_z (zD)P(z)- M_1(z) \nonumber \\
&=& \lambda_z \nu_z({\M}_z^{(1)}h_z) - M_1(z)\sim (\lambda_z-1)\, M_1(z),
\end{eqnarray}
so that $\lim_{z\uparrow 1}K(z)/(1-z)=0$ for all $d>-1$ and $\lim_{z\uparrow 1}K(z)/(1-z)^2<\infty$ for $d>0$.
Hence, for all $d>-1$, we have $\lim_{z\uparrow 1}J(z)=0$.  
This gives the asymptotic behaviour of $Q(z)$
and the nature of the singularity at $z=1$ follows from Proposition \ref{generfcts}.

\noindent
More generally, if $m-1<d\leq m$ for some $m>0$ then
both $(zD)^{m+1}\lambda_z$ and $M_{m+1}(z)$ diverge as $z\uparrow 1$  but 
again from Theorem \ref{analytic} we have
\be\label{kemme}
(zD)^{m+1}K(z)=(zD)^{m+1}\lambda_z - M_{m+1}(z)={\cal O}(1).
\ee
This gives
$(zD)^{m}J(z) = o\, (1)$ and $(zD)^{m+1}J(z)={\cal O}(M_{m+1}(z))$  as $z\uparrow 1$, so that the behaviour of $U(z)$
follows from the identities
\be\label{dec2}
U(z)=Q(z)\left( M_1- {1-e^{P(z)}\over 1-z}\right) =B(z) + M_1A(z) J(z)
\,.
\ee
More specifically,
differentiating the term $A(z) J(z)$ we get
\be\label{indu}
(zD)\left[ A(z) J(z)\right]=A(z) \cdot\biggl\{\left( (zD)\,e^{P(z)}\right) A(z) J(z)+(zD)J(z)\biggr\}
\ee
which, by the above, is ${\cal O}\left(M_{2}(z)/M_1(z)\right)$  as $z\uparrow 1$.
This implies that if $-1<d\leq 0$ the term $A(z) J(z)$ is continuous at $z=1$ and 
$(zD)\left[ A(z) J(z)\right] = {\cal O}\left( (1-z)^{-1}\right)$ when $z\uparrow 1$.
For $0<d\leq 1$ we have $A(z) J(z) = o\, (1)$  and 
$(zD)\left[ A(z)J(z)\right] ={\cal O}\left(M_2(z)\right)=o\, ((1-z)^{-1})$ when $z\uparrow 1$. In particular, the
last expression is integrable in a neighbourhood of $z=1$. More generally, setting
$H_\ell(z):=A(z)^{-1}(zD)^\ell\left[ A(z) J(z)\right]$, so that
$H_{\ell+1}(z)=(zD)H_\ell (z) + H_\ell (z) \, (zD)\left[\log A(z)\right]$,
one readily sees inductively that if $m-1<d\leq m$ for some $m>0$ then the term
$A(z)\cdot H_m(z)$ is ${\cal O}((zD)^{m+1}J(z))={\cal O}\left(M_{m+1}(z)\right)$ as  $z\uparrow 1$, 
and thus integrable as above. $\qed$

\section{Renewal vs scaling and mixing properties}\label{rensca}
We are now going to make use of the results of the previous sections to study 
the generating function of the numbers $u_n$'s (see (\ref{u_n})), that is the function
\be
S(z):=\sum_{n=0}^{\infty} z^n \, u_n.
\ee
We have the following
\begin{proposition}\label{genthm} $S(z)$ is holomorphic in $\D$ and extends continuously to
$\Dc \setminus \{1\}$. Moreover,
for $z$ in some neighbourhood of $z=1$, with $z\ne 1$, it can be written as
$$
S(z)=\nu ({\P}_z h)\cdot Q(z)+ \nu \left( (1-{\Nc}_z)^{-1}h\right),
$$
where ${\cal P}_z = h_z\cdot \nu_z$ is the spectral projection of ${\M}_z$ corresponding
to the eigenvalue $\lambda_z$.
\end{proposition}
{\sl Proof.}
We start with the identity
\be\label{formalident0}
S(z) = \sum_{n=0}^{\infty} z^n \, \nu (1_{A_1}\cdot {\L}^n (1_{A_1}\cdot e)\, ) 
= \nu \left( 1_{A_1}\cdot (1-z{\L})^{-1} (1_{A_1}\cdot e)\, \right),
\ee
so that Corollary \ref{decco} yields the first part of the theorem along with the expression, 
valid in some neighbourhood of $z=1$,
\be\label{e1}
S(z)=  \nu (1_{A_1}\cdot e_z)\cdot \nu_z(1_{A_1}\cdot e)\cdot Q(z)  +  {R}(z)
\ee
with  
${R}(z)=\nu \biggl( 1_{A_1}\cdot (1-z{\L}_0)^{-1}(1-{\Nc}_z)^{-1}(1_{A_1}\cdot e)\, \biggr)$.
Now notice that
\begin{eqnarray}
\nu (1_{A_1}\cdot e_z)&=&\nu ({\M}(1_{A_1}\cdot e_z))=\nu ({\L}_1 e_z) 
=\nu({\L}_1(1-z{\L}_0)^{-1}h_z)\nonumber \\
&=&\nu\left({{\M}_z h_z\over z}\right)={\lambda_z \over z}\,\nu (h_z), 
\end{eqnarray}
and also
\be
\nu_z (1_{A_1}\cdot e)=\nu_z \left({{\M}_z(1_{A_1}\cdot e)\over \lambda_z}\right)=
{z\over \lambda_z}\, \nu_z ({\L}_1 e) ={z\over \lambda_z}\, \nu_z(h).
\ee
Therefore $\nu (1_{A_1}\cdot e_z)\cdot \nu_z(1_{A_1}\cdot e) = \nu (h_z)\cdot \nu_z(h)$.
Using (\ref{anarel}) we then have
\be
\nu (h_z)\cdot \nu_z(h) =\nu ({\cal U}_zh)\cdot \nu ({\cal U}_z^{-1}h)=
\nu \left({\cal U}_z \, h \, \nu \,{\cal U}_z^{-1}\, h\right)=\nu({\cal P}_z h).
\ee
Finally, since ${\M}_z$ and ${\Nc}_z$ are commuting we have
\begin{eqnarray}
R(z)&=&
\nu \biggl( 1_{A_1}\cdot (1-z{\L}_0)^{-1}(1-{\Nc}_z)^{-1}(1_{A_1}\cdot e)\, \biggr)\nonumber\\
&=&
\nu \biggl( {\L}_1\, (1-z{\L}_0)^{-1}(1-{\Nc}_z)^{-1}(1_{A_1}\cdot e)\, \biggr)\nonumber\\
&=&
\nu \biggl( z^{-1}{\M}_z(1-{\Nc}_z)^{-1}(1_{A_1}\cdot e)\, \biggr)\nonumber\\
&=&
\nu \biggl( (1-{\Nc}_z)^{-1}z^{-1}{\M}_z(1_{A_1}\cdot e)\, \biggr)\nonumber\\
&=&
\nu \biggl( (1-{\Nc}_z)^{-1} h\, \biggr).\qquad\qquad\qquad\qquad \qed\nonumber
\end{eqnarray}
The above and formula (\ref{dec1}) allow us to write 
\be\label{augh}
S(z)=A(z)+C(z),
\ee
where $A(z)$ is defined in (\ref{generatingfct1}) and $C(z)$ 
is a function which is holomorphic in $\D$ and in a neighbourhood of $z=1$ can be written as
\be\label{difference}
C(z)=A(z)\cdot J(z) +Q(z)\cdot \Delta (z)  +\nu \biggl( (1-{\Nc}_z)^{-1} h\, \biggr)
\ee
where we have set $\Delta (z):=\nu ({\P}_z h)-1$. We now investigate the behaviour of 
$C(z)$ on the unit circle.
To this end, we can evaluate the function $\Delta (z)$ when $z\uparrow 1$ as follows. 
First, using (\ref{estim}), (\ref{asymp}) 
along with $ h\asymp 1$, and Corollary \ref{singular} we get
\be
 \Vert {\cal M}_z-{\cal M}_1\Vert_\t \leq \sum_{n\geq 1}(1-z^n)\Vert {\Mc}^{(n)}\Vert_\t
\leq C\, (1-\lambda_z).
\ee
On the other hand we have that $|\Delta (z)|Ê\leq \Vert {\P}_z-{\cal P}_1\Vert_\t$.
By the spectral properties of ${\M}_z$ the last quantity is of the same order as
$\Vert {\cal M}_z-{\cal M}_1\Vert_\t$ and thus, by the above, of order not larger than $1-\lambda_z=Q(z)^{-1}$. 
Furthermore, since ${\Nc}_1 h=0$
the same holds true for the quantity $\nu \biggl( (1-{\Nc}_z)^{-1} h\, \biggr)$.
Now notice that by (\ref{momenti}) we have, for $z\uparrow 1$,
\be\label{deltaprime}
\Delta(z) =\nu({\cal P}_zh)-1\sim
\nu({\cal M}_zh-{\cal M}_1h)\sim \nu({\cal M}_z^{(1)}h)(z-1) = M_1(z)(z-1)
\ee
and therefore $Q(z)\cdot \Delta(z)$ is continuous at $z=1$ with
\be\label{-1}
\lim_{z\uparrow 1} Q(z)\cdot \Delta(z) =-1.
\ee
In particular, if $d>0$ then $\Delta (z)$ is differentiable at $z=1$ and $(zD)\Delta(z)|_{z=1}=M_1$.
Differentiation of $Q(z)\cdot \Delta(z)$ leads to the expression
\be\label{der}
(zD)\left[ Q(z)\cdot \Delta(z)\,\right] =Q(z) \cdot\biggl\{\left( (zD)(e^{P(z)})\right) Q(z)\Delta (z)
+(zD)\Delta(z)\biggr\},
\ee
so that proceeding as in the proof of Corollary \ref{singular} we obtain that for all $d>-1$ the function
$(zD)[ Q(z)\cdot \Delta(z)]$ is ${\cal O}\left( M_2(z)/M_1(z)\right)$ as $z\uparrow 1$. 
A similar reasoning shows that for all $d>-1$ the function
$(zD)\nu \biggl( (1-{\Nc}_z)^{-1} h\, \biggr)$ is ${\cal O}\left( M_1(z)\right)$ as $z\uparrow 1$. 

\noindent
A first consequence of the above discussion is that the function
$C(z)$ is holomorphic in $\D$ and extends continuously to
$\Dc$. In particular, the function $C(e^{i\phi})$ is integrable
on $[-\pi,\pi]$. By the theorem that Fourier coefficients 
tend to zero (see \cite{Zig}, p.45) we then have that the coefficient of $z^n$ 
in the power series expansion of $C(z)$ is $o\, (1)$.
\vskip 0.2cm
\noindent
Let us consider first the infinite measure case, i.e. $-1<d\leq 0$.
From the previous discussion and Proposition \ref{generfcts} we have that the coefficient of $z^n$ in the power series expansion of $S(z)$ 
tends to zero as $n\to \infty$. But we can say more. Multiplying (\ref{augh}) by the function $D(z)$ dealt with in Proposition \ref{generfcts} 
we obtain $D(z)\cdot S(z)= (1-z)^{-1}+D(z)\cdot C(z)$. This function is continuous on $\Dc\setminus\{1\}$. Moreover, from Corollary \ref{singular} and the above 
it follows that $D(z)\cdot C(z) = {\cal O}\left(M_2(z)(1-z)\right)$ in a neighbourhood of $z=1$ so that
$\lim_{r\uparrow 1} \int_{-\delta}^\delta |D(re^{i\theta})\cdot C(re^{i\theta})| d\theta = {\cal O}(1)$ 
for all $\delta >0$. Therefore the function $D(z)\cdot S(z)-(1-z)^{-1}$ is (absolutely) integrable
on the unit circle $|z|=1$, and the coefficients
of its power series expansion tend to zero. Whence
$\sum_{k=0}^n d_k\, u_{n-k} = 1 + o\, (1)$ and, by (\ref{generatingfct1}),
$\sum_{k=0}^n d_k\, (u_{n-k}-a_{n-k}) = o\, (1)$. By a simple lemma (see \cite{Chu}, Chap. I.5, Lemma A)
and  Proposition \ref{generfcts} we then obtain that for $-1<d\leq 0$,
\be
u_n \sim a_n \approx n^d,\qquad n\to \infty,
\ee
where, again, for $d=0$ this relation means that $u_n$ tends to zero slower than any inverse power of $n$.

\noindent
We now examine the finite measure case $d>0$. Here, putting together (\ref{augh}), Proposition \ref{generfcts} and the vanishing of
the coefficients of the power series expansion of $C(z)$ established so far we get
\be \label{renthm}
u_n \to {1\over M_1}\quad\hbox{as}\quad n\to \infty.
\ee
In view of (\ref{passages}), (\ref{renthm}) can be
regarded as a {\sl renewal property} as well as a {\sl mixing property} for the dynamical system
$(\O,T,\mu)$. Indeed, in this case we can rewrite
(\ref{renthm}) as 
\be\label{mix}
{\hat \mu} (A_1\cap T^{-n}A_1) \to ({\hat \mu}(A_1))^2\quad\hbox{as}\quad n\to \infty, 
\ee
where ${\hat \mu}$ is the probability measure $\mu/M_1$. 

\noindent
In order to investigate the speed of convergence in the limit (\ref{mix}) we shall consider the numbers 
\be
v_n=M_1u_n-1={
{\hat \mu} (A_1\cap T^{-n}A_1) - ({\hat \mu}(A_1))^2\over({\hat \mu}(A_1))^2}
\ee
along with their
generating function $S_1(z)=\sum_{n=0}^{\infty} z^n \,v_n$. 
By Proposition \ref{genthm}, $S_1(z)$ is holomorphic in $\D$, extends continuously to
$\Dc \setminus \{1\}$ and can be written as
\be\label{S_1}
S_1(z)= B(z)+ M_1\, C(z),
\ee
where $B(z)$ is defined in (\ref{generatingfct2}) 
and its behaviour is described in 
Proposition \ref{generfcts} while $C(z)$, defined in (\ref{augh}), 
is continuous on the unit circle {\sl and} its derivative is integrable on the same domain. 
Therefore by the same reasoning as above the coefficient of $z^n$ 
in the power series expansion of $C(z)$ is $o\, (n^{-1})$. Comparing to Proposition \ref{generfcts} we then see
that for all $d$ satisfying $0<d\leq 1$ one has $v_n \sim b_n$.
More generally,
using 
(\ref{derivatives}), (\ref{formalident0}) and (\ref{augh}) we see that if $d>m-1$ for some $m>0$ then
$(zD)^mS(z)$ still extends continuously to $\Dc\setminus\{1\}$ 
and clearly the same can be said for $(zD)^m B(z)$. 
One can now proceed inductively, as in the proof
of Corollary \ref{singular},
and realize that if $d>m-1$ for some $m>0$ then the function
$(zD)^mC(z)$ is integrable on the unit circle so that the coefficient of $z^n$ in its power series expansion is 
$o(1)$ and thus that in the power series expansion of $C(z)$ is $o(n^{-m})$. Again comparing to Proposition \ref{generfcts} one then sees
that for all $d$ satisfying $m-1<d\leq m$ one has $v_n \sim b_n$.

\noindent
We have thus proved that the behaviour of both $u_n$ and $v_n$ for $(\O,T,\mu)$ is the same 
as that of the corresponding quantities $a_n$ and $b_n$ for the Markov approximation:

\begin{theorem} \label{rateofmixing} For $-1<d\leq 0$ we have
$$
u_n=\mu(A_1\cap T^{-n}A_1) \sim a_n \approx n^d
$$
wheras for $d>0$ we have
$$
v_n=M_1u_n-1={
{\hat \mu} (A_1\cap T^{-n}A_1) - ({\hat \mu}(A_1))^2\over({\hat \mu}(A_1))^2}\sim b_n \sim {\hat \mu} (\tau > n) \approx n^{-d}.
$$
\end{theorem}
If we now consider the partition set $A_\ell$ for some $\ell> 1$, we have that 
$\mu(A_\ell\cap T^{-n}A_\ell)=0$ for $0<n<\ell$ and a computation similar to that 
in the proof of Proposition \ref{genthm} leads to the following identity 
(valid for $z$ the vicinity of $1$, $z\ne 1$):
\be
\sum_{n\geq \ell} z^n \, \mu(A_\ell\cap T^{-n}A_\ell) =
\nu \left( {\Mc}_{z, \ell}\, {\cal P}_z\, h_\ell\right)\cdot Q(z)
+\nu({\Mc}_{z,\ell}\, (1-{\Nc}_z)^{-1} h_\ell)
\ee
where 
\be
{\Mc}_{z, \ell}\, u (\o) :=\sum_{k=\ell}^{\infty} z^{k}\, 
\left({\cal L}_1 {\cal L}_0^{k-1}u\right) (\o)\quad\hbox{and}\quad h_\ell :={\Mc}_{1,\ell}\, h.
\ee
Notice that the terms $\nu \left( {\Mc}_{z, \ell}\, {\cal P}_z\, h_\ell\right)$ and
$\nu({\Mc}_{z,\ell}\, (1-{\Nc}_z)^{-1} h_\ell)$
are ${\cal O}(z^\ell)$ when $z\downarrow 0$ and satisfy  
\be\label{lim1}
\lim_{z\uparrow 1} \nu \left( {\Mc}_{z, \ell}\, {\cal P}_z\, h_\ell\right) =
\left( \nu(h_\ell)\right)^2=\left( \sum_{k\geq \ell}\rho(A_k)\right)^2 =
\left(\mu(A_\ell)\right)^2,
\ee
and, since ${\Nc}h_\ell = (1-\mu (A_\ell) )\, h$, 
\be\label{lim2}
\lim_{z\uparrow 1} \nu({\Mc}_{z,\ell}\, (1-{\Nc}_z)^{-1} h_\ell)=\mu(A_\ell)(1-\mu(A_\ell) ),
\ee
respectively.

\noindent
On the other hand, a short calculation yields the decompositions
\be\label{genscal}
\sum_{n\geq \ell} z^n \, {\mu(A_\ell\cap T^{-n}A_\ell)\over (\mu(A_\ell))^2}  =  z^\ell \cdot A(z) +  
{1\over ({\hat \mu}(A_\ell))^2}\,C_\ell (z),
\ee
and, for $d>0$,
\be\label{gencorr}
\sum_{n\geq \ell} z^n \,\left( {
{\hat \mu}(A_\ell\cap T^{-n}A_\ell) -({\hat \mu}(A_\ell))^2\over ({\hat \mu}(A_\ell))^2 } \right)
=z^\ell \cdot B(z) + {1\over M_1({\hat \mu}(A_\ell))^2}\, C_\ell (z),
\ee
where $C_\ell(z)$ is holomorphic in $\D$, extends continuously to
$\Dc \setminus \{1\}$, and
for $z$ in a neighbourhood of $z=1$ gets the expression
\begin{eqnarray}\label{Cell}
C_\ell (z) &=& z^\ell \cdot({\hat \mu}(A_\ell))^2\, A(z) J(z)+ 
Q(z)\cdot 
\left[\nu \left( {\Mc}_{z, \ell}\, {\cal P}_z\,h_\ell\right)- z^\ell \,(\mu(A_\ell))^2
\right]\nonumber\\ & &\;\;\;\, +\; \; \nu({\Mc}_{z,\ell}\, (1-{\Nc}_z)^{-1}h_\ell)\; .
\end{eqnarray}
Therefore, the same reasoning leading to Theorem \ref{rateofmixing} can be applied here
to give
\be\label{x1}
{\mu(A_\ell\cap T^{-n}A_\ell)\over (\mu(A_\ell))^2} \sim u_{n-\ell} \sim u_n,
\ee
and, for $d>0$, 
\be\label{x2}
{{\hat \mu} (A_\ell\cap T^{-n}A_\ell) - ({\hat \mu}(A_\ell))^2\over ({\hat \mu}(A_\ell))^2} \;\sim \;
v_{n-\ell} \sim v_n,
\ee
where the last asymptotic equivalences in both expressions
hold for each fixed $\ell\in \Z_+$. 
It is now an easy matter to realize that for any given $\ell \geq 1$ and any Borel 
set $E\subseteq A_\ell$ with $\mu(E)>0$
one has
\be
\sum_{n\geq \ell} z^n \, {\mu(E\cap T^{-n}E)\over (\mu(E))^2}  =  z^\ell \cdot A(z) +  
{1\over ({\hat \mu}(E))^2}\, C_E (z),
\ee
and, for $d>0$,
\be
\sum_{n\geq \ell} z^n \,\left( {
{\hat \mu}(E\cap T^{-n}E) -({\hat \mu}(E))^2 \over ({\hat \mu}(E))^2} \right)
=z^\ell \cdot B(z) + {1\over M_1 \,({\hat \mu}(E))^2}\, C_E (z),
\ee
where the expression of $C_E(z)$ in a neighbourhood of $z=1$ is as in (\ref{Cell}) 
with $\mu(A_\ell)$ replaced by $\mu(E)$ and the quantities ${\Mc}_{z,\ell}$, $h_\ell$
replaced by ${\Mc}_{z,E}$, $h_E$ given by
\be
{\Mc}_{z, E}\, u (\o) :=\sum_{k=\ell}^{\infty} z^{k}\, 
\left({\cal L}_1 {\cal L}_0^{k-1}(u\cdot 1_{T^{-(k-\ell)}E\cap A_k})\right) (\o)\quad\hbox{and}\quad
 h_E :={\Mc}_{1,E}\, h.
\ee
Moreover, from these relations we see that the limits corresponding to (\ref{lim1}) and (\ref{lim2}) hold with 
$\mu(A_\ell)$ replaced by $\mu(E)$.
The uniform distortion property (\ref{asymp}) now allows us to repeat exactly the same reasoning as above
to obtain, for any fixed $\ell \in \Z_+$ and $E\subseteq A_\ell$ with $\mu (E)>0$,
\be\label{decay1}
{\mu(E\cap T^{-n}E)\over (\mu(E))^2} \sim u_{n-\ell} \sim u_n,
\ee
and, for $d>0$, 
\be\label{decay2}
{{\hat \mu} (E\cap T^{-n}E) - ({\hat \mu}(E))^2\over ({\hat \mu}(E))^2} \;\sim \;
v_{n-\ell} \sim v_n.
\ee
In an entirely analogous way one shows that (\ref{decay1}) and (\ref{decay2}) hold true
for $E\subseteq \cup_{\ell\in J}A_\ell$
where $J\subset \Z_+$ is any given finite set.

\noindent
Let ${\cal B}(\O)$ be the Borel $\sigma$-algebra on $\O$.
Given $E\subset {\cal B}(\O)$, for $-1<d\leq 0$ we define the {\sl scaling rate} $\s_n(E)$ of $E$ as
\be
\s_n(E) := {{ \mu} (E\cap T^{-n}E)\over ({\mu}(E))^2},
\ee
and, for $d>0$,
the {\sl mixing rate} $\mu_n(E)$ of $E$ as
\be
\mu_n(E) := {{\hat \mu} (E\cap T^{-n}E) - ({\hat \mu}(E))^2\over ({\hat \mu}(E))^2}\cdot
\ee
This quantities
are not uniform in $E\subset {\cal B}(\O)$. As already noted, the last asymptotic equivalences
in (\ref{decay1}) and (\ref{decay2}) follow for each fixed $\ell\in \Z_+$, but not uniformly in $\ell$.
To recover uniformity, we consider the set
$D_N=\{\o \in \O : \tau (\o) > N\}= \cup_{\ell >N}A_\ell$ and define
\be
B_{+}:= \cup_{N} \{E\in {\cal B}(\O): \, \mu (E)>0, 
\, E \subseteq \O \setminus  D_N\,\}.
\ee
An easy consequence of the above discussion is the following
\begin{lemma}  Let $E,F\in B_+$. If $-1<d\leq 0$ then $\s_n(E) \sim \s_n(F)$. If instead $d>0$ then $\mu_n(E) \sim \mu_n(F)$.
\end{lemma}
\noindent
Therefore, one may give the following definition,
\begin{definition} For $-1<d\leq 0$ the scaling rate $\s_n(T)$ of $(\O,T,\mu)$ is the rate of 
asymptotic decay of the sequences $\{\s_n(E)\}$, with $E\in B_+$. Similarly, for $d>0$, 
the mixing rate $\mu_n(T)$ of $(\O,T,\mu)$
is the rate of asymptotic decay of the sequences $\{\mu_n(E)\}$, with $E\in B_+$.
\end{definition}
We summarize the previous findings in the following
\begin{theorem}\label{mr}
For $-1<d\leq 0$ we have
$\s_n(T) \approx n^d$, 
wheras for $d>0$ we have
$\mu_n (T)= {\hat \mu} (\tau > n) \approx n^{-d}$.
\end{theorem}
\begin{remark}\label{8} {\rm
Let $E\subset {\cal B}$, with $0<m(E)<1$, and $E^c=\O\setminus E$. It is an easy 
observation that (see \cite{Is1})
$$
{\hat \mu}(E\cap T^{-n}E)-({\hat \mu}(E))^2={\hat \mu}(E^c\cap T^{-n}E^c)-({\hat \mu}(E^c))^2.
$$
Now assume that $E\subset B_+$ so that $E^c \not \subset B_+$.
The above identity and Theorem \ref{mr} then imply that 
$$
\mu_n (E^c)\sim \left( {{\hat \mu}(E)\over {\hat \mu}(E^c)}\right)^2\, {\hat \mu} (\tau > n),
\quad E\subset B_+,\quad 0<m(E)<1.
$$
This can be used, for instance, to evaluate the mixing rate of the set $D_N=\cup_{\ell>N}A_\ell=\{\tau >N\}$, for
any $N\in \Z_+$.}
\end{remark}
\begin{remark}
{\rm When $M_1=\infty$ the scaling rate $\s_n(T)$ satisfies 
\be
\sum_{k=1}^{n}\s_k(T) \cdot \sum_{k=1}^{n}\mu(\tau=k) 
\sim n\quad\hbox{as}\quad n\to \infty.
\ee
This establishes a (asymptotic) relation between the scaling rate
and the behaviour of the partial sums $\sum_{k=1}^{n}\mu(\tau=k)$
which, in turn,
is related to quantities such as the {\sl wandering rate} and the {\sl return sequence} that naturally arise in the context of 
infinite ergodic theory and for which we refer to \cite{Aa} (see also \cite{Is1}).}
\end{remark}
\begin{remark} {\rm The part for $d>0$ of
 Theorem \ref{mr} was proved in \cite{Is1} for the Markov chain 
of Example 1 (see also [FL]).  We point out that it gives 
the (asymptotically) exact rate of mixing for $(\O,T,\mu)$, not just a bound for it, and
can be viewed as a statement about the decay of correlations
for test functions as simple as indicators of sets in $B_+$. 
This makes the mixing rate (as defined above) determined by nothing else than the
distribution of return times: $ {\hat \mu} (\tau > n)$.
On the other hand, when dealing with correlation functions of
a broader class of observables, one expects a richer behaviour
depending also of the smoothness properties of the functions involved along with their behaviour in the vicinity of
the fixed point $\{0^\infty\}$.
In particular one may obtain faster decays.}
\end{remark}

\section{Weak Bernoulli and polynomial cluster properties}\label{poldecay}

In this Section we shall consider only the ergodic case $d>0$ and obtain some further properties of the probability measure 
${\hat \mu}$. To this end we first observe that
the reasoning of the previous Section can be easily extended to 
study the behaviour of the quantity 
$\mu(A_\ell\cap T^{-n}A_k)$,  with $k\ne \ell$, yielding (notation as in the previous Section):
\be
\sum_{n\geq \ell} z^n \, \mu(A_\ell\cap T^{-n}A_k) = z^{\ell-k}\left[
\nu \left( {\Mc}_{z,k}\, {\cal P}_z\, h_\ell\right)\cdot Q(z)
+\nu({\Mc}_{z,k}\, (1-{\Nc}_z)^{-1} h_\ell)\right]
\ee
and therefore
\be\label{gencorr}
\sum_{n\geq \ell} z^n \,\left( {
{\hat \mu}(A_\ell\cap T^{-n}A_k) -{\hat \mu}(A_k){\hat \mu}(A_\ell)\over {\hat \mu}(A_k){\hat \mu}(A_\ell) } \right)
=z^\ell \cdot B(z) + {1\over M_1{\hat \mu}(A_k){\hat \mu}(A_\ell)}\, C_{k,\ell} (z),
\ee
where $C_{k,\ell}(z)$ is holomorphic in $\D$, extends continuously to
$\Dc \setminus \{1\}$, and
for $z$ in a neighbourhood of $z=1$ writes
\begin{eqnarray}\label{Ckell}
C_{k,\ell} (z) &=& z^\ell \cdot {\hat \mu}(A_k){\hat \mu}(A_\ell)\, A(z) J(z)+ 
Q(z)\cdot 
\left[z^{\ell-k}\nu \left( {\Mc}_{z, k}\, {\cal P}_z\,h_\ell\right)- z^\ell \, \mu(A_k)\mu(A_\ell)
\right]\nonumber\\ & &\;\;\;\, +\; \; z^{\ell-k}\cdot \nu({\Mc}_{z,k}\, (1-{\Nc}_z)^{-1}h_\ell)\; .
\end{eqnarray}
The same reasoning as above then yields
\be
{{\hat \mu} (A_\ell\cap T^{-n}A_k) - {\hat \mu}(A_k){\hat \mu}(A_\ell)\over {\hat \mu}(A_k){\hat \mu}(A_\ell) }\;\sim \;
v_{n-\ell}\; .
\ee
Now recall that the partition sets $A_\ell$ are particular cylinder sets: $A_\ell =[0^{\ell-1}1]$. On the other
hand, given an arbitrary $\ell$-cylinder set $E=[\o_0,\cdots ,\o_{\ell-1}]$ with $E\ne D_\ell\equiv [0,\cdots , 0]=\{\tau >\ell\}$, 
then either $E=A_\ell$ or
$E\subset A_i$ for some $i<\ell$. 
Thus, a straightforward consequence of the previous discussion is that 
we can find a positive constant 
$C$ such that for any pair of
$\ell$-cylinders $E,F$ both $\not=D_\ell$ we have
\be \label{basineq}
|{\hat \mu} (E\cap T^{-n}F) - {\hat \mu}(E){\hat \mu}(F)| \leq C\, {\hat \mu}(E){\hat \mu}(F)\, v_{n-\ell}\, .
\ee
Moreover, noting that for all $F\subset B_+$ (see also remark \ref{8})
\begin{eqnarray}
|{\hat \mu} (D_\ell\cap T^{-n}F) - {\hat \mu}(D_\ell){\hat \mu}(F)|&=&|{\hat \mu} (D_\ell^c\cap T^{-n}F) -
 {\hat \mu}(D_\ell^c){\hat \mu}(F)|\nonumber \\
&\leq& C\, {\hat \mu}(D_\ell^c){\hat \mu}(F)\, v_{n-\ell},
\end{eqnarray}
(the same inequality holds with $D_\ell$ and $F$ interchanged) and
summing over all  $\ell$-cylinder sets $E$, $F$ we get (recall that ${\hat \mu}(\O)=1$)
\be\label{wber}
\sum_{E,F} |{\hat \mu} (E\cap T^{-n}F) - {\hat \mu}(E){\hat \mu}(F)| \leq 
C(1+2\,{\hat \mu}(\{\tau \leq \ell\}))\cdot v_{n-\ell}
\ee
Therefore, since $v_n\approx n^{-d} \to 0$, we obtain (see \cite{Bo} for definitions)
\begin{theorem}\label{wb}
The partition $\{[0],[1]\}$ is weak-Bernoulli (and hence Bernoulli) for $T$ and ${\hat \mu}$.
\end{theorem}
We now proceed as in \cite{Bo} to carry the polynomial convergence (\ref{basineq}) over to functions
of polynomial variation. Given $a >0$ let ${\cal H}_a(\O)$ be the family of $f\in C(\O)$ with
$\var_n f \leq C n^{-a}$. The space ${\cal H}_a$ becomes a Banach space under the norm
\be 
\Vert f \Vert_a = |f|_\infty + \sup_{n >0}\left\{ n^a \cdot \var_n f \right\}.
\ee
We have the following polynomial cluster property:
\begin{theorem}\label{pcp} For $d>0$ and fixed $a\geq  d$ there is a constant $D>0$ such that
$$
|{\hat \mu} (f\cdot g\circ T^{n}) - {\hat \mu}(f){\hat \mu}(g)|\leq D \Vert f\Vert_a \Vert g\Vert_a n^{-(d-\epsilon)}
$$
for all $n>0$, $\epsilon >0$ and $f,g \in {\cal H}_a$.
\end{theorem}
{\sl Proof.} 
Given a pair $f,g \in {\cal H}_a(\O)$ we let $f_\ell$ and $g_\ell$ be the conditional expectations of 
$f$ and $g$ wrt the $\s$-algebra generated by $\ell$-cylinder sets, i.e.
\be f_\ell = \sum_E f_E \, 1_E \quad\hbox{with}\quad f_E={1\over {\hat \mu}(E)}\int_E f(\o) {\hat \mu}(d\o)
\ee
and similar expression for $g$. We have ${\hat \mu}(f_\ell)={\hat \mu}(f)$ and 
$|f-f_\ell| \leq \Vert f\Vert_a \cdot n^{-a}$. A routine estimate (see \cite{Bo}, p.39) then gives
\be
|{\hat \mu} (f\cdot g\circ T^{n}) - {\hat \mu}(f){\hat \mu}(g)| \leq 
|{\hat \mu} (f_\ell\cdot g_\ell\circ T^{n}) - {\hat \mu}(f_\ell){\hat \mu}(g_\ell)| +
 2 n^{-a} \Vert f\Vert_a \Vert g\Vert_a.
\ee
Moreover, by (\ref{wber}), we have
\begin{eqnarray}
|{\hat \mu} (f_\ell\cdot g_\ell\circ T^{n}) - {\hat \mu}(f_\ell){\hat \mu}(g_\ell)| &\leq &
\sum_{E,F} |f_E| \, |g_F|\, |{\hat \mu} (E\cap T^{-n}F) - {\hat \mu}(E){\hat \mu}(F)|\nonumber \\
&\leq & C\, |f|_\infty\, |g|_\infty\, (1+2\,{\hat \mu}(\{\tau \leq \ell\}))\cdot v_{n-\ell}\nonumber \\
&\leq & 3\, C\, \Vert f\Vert_a \Vert g\Vert_a \,  v_{n-\ell}\, . \nonumber 
\end{eqnarray}
The assertion now follows putting together the last two inequalities with the choice $\ell =[n/2]$. $\qed$

\section{Zeta functions}\label{zetafunc}
We consider the dynamical zeta functions $\zeta (\p,z)$ and $\zeta (\psi,w)$
associated to the pairs $(T,\varphi)$ and $(S,\psi)$, respectively 
(here, as in Section \ref{opvalpowser}, $S$ denotes the induced map $T^\tau$ when acting on the symbol space $\S$)
and defined by 
the following formal series (see \cite{Ru1}, \cite{Ru2}, \cite{Ba2}):
$$
\zeta (\p, z) = \exp \sum_{n=1}^{\infty} {z^n\over n} Q_n \quad\hbox{and}\quad
\zeta (\psi,w) = \exp \sum_{n=1}^{\infty} {w^n\over n} Z_n,
$$
where the `partition functions' $Q_n$ and $Z_n$ are given by
$$
Q_n =\sum_{\o=T^n\o} \prod_{k=0}^{n-1}\p(T^k\o)
\quad\hbox{and}\quad
Z_n =\sum_{\s=S^n\s} \prod_{k=0}^{n-1}\psi(S^k\s).
$$
\begin{remark}
{\rm Notice that
the correspondence $\s_j(\o)$ defined in (\ref{passages0}) 
may associate to a periodic sequence an eventually periodic one.
More precisely, 
let $\o \in \O_0$ be a periodic sequence of period $n$ for the shift
$T$. We write it in the form 
$\o=({ \overline {\o_0 \o_1 \dots \o_{n-1} }})$. Now, 
if $\o_0=0$, 
it may happen that, for some $k\geq 1$,
$$
\o_0 \o_1 \dots \o_{n-1}\,  = \,{\underbrace{00\dots 01}_{l_0\geq 2}}\,\,
{\underbrace{1\dots 1}_{r_1\geq 0}}
\,\, {\underbrace{00\dots 01}_{l_1\geq 2}}\,\,\ldots \,\,
{\underbrace{1\dots 1}_{r_k\geq 0}}\,\,
{\underbrace{00\dots 0}_{l_k\geq 1}}
$$
and one would find the eventually periodic
sequence $\s(\o)=({\s_0 \overline {\s_1 \s_2 \dots \s_{m} }})$ where
$$
\s_0 = l_0\quad\hbox{and}\quad 
\s_1 \s_2 \dots \s_{m}\, = \, {\underbrace{1\dots 1}_{r_1\geq 0}}
\,\, l_1\,\,\ldots \,\,
{\underbrace{1\dots 1}_{r_k\geq 0}}\,\, (l_k+ l_0)
$$
whose ultimate period is $m= k+ r_1+\dots +r_k$ and which satisfies 
$n=\s_1+\dots +\s_m$.
Thus, in the latter case, in order to obtain a correspondence 
between periodic
sequences it is necessary to apply the 
aformentioned rule to some iterate
of the original sequence. }
\end{remark}
Let us now examine how $\zeta (\p,z)$ and $\zeta (\psi,w)$ are related
to one another. 
To this end we observe that if $\o\in \O_0$ is a
periodic sequence of period $n$ for the
shift $T$ and $\s (\o )\in \S$ is the periodic sequence of period 
$m=\sum_{k=0}^{n-1}1_{A_1}(T^k\o)$ 
corresponding to $\o$ or to some iterate of it 
(see the above remark), then we have $n=\sum_{j=0}^{m-1}\s_j$ and
$$
\prod_{k=0}^{n-1}\p(T^k\o) = \prod_{k=0}^{m-1}\psi(S^k\s).
$$
Using this facts we write $Q_n$ as follows:
$$  
Q_n = 1+ \sum_{m=1}^n {n\over m}
\sum_{\scriptstyle \s=S^m\s}
\prod_{k=0}^{m-1}\psi(S^k\s)
$$  
where the $1$ comes from the fixed point $0^{\infty}$.
The second sum ranges over the $n-1 \choose m-1$ ways to write the integer $n$
as a sum of $m$ positive integers, counting all permutations.
Therefore, 
\begin{eqnarray}
\sum_{n=1}^{\infty}{z^n\over n} Q_n &=& \log ({1\over 1- z}) +
\sum_{n=1}^{\infty}\sum_{m=1}^{n}
{1\over m} \sum_{\scriptstyle \s=S^m\s} 
z^n \prod_{k=0}^{m-1}\psi(S^k\s) \nonumber \\
&=& \log ({1\over 1- z}) +
\sum_{\ell=1}^{\infty} {1\over \ell} \sum_{\scriptstyle \s=S^\ell\s }
 z^{\sum_{j=0}^{\ell-1}\s_j}\prod_{k=0}^{\ell-1}\psi(S^k\s) \nonumber
\end{eqnarray}
Putting together these observations we have the following result, which
can be viewed as the counterpart of Proposition \ref{identity}:
\begin{proposition} \label{twozeta} Consider the `grand partition function' $\Xi_\ell(z)$ defined by
$$
\Xi_\ell(z) := \sum_{\scriptstyle \s=S^\ell\s }
z^{\sum_{j=0}^{\ell-1}\s_j}\prod_{k=0}^{\ell-1}\psi(S^k\s)
$$
and
the two-variable
zeta function given by
$$
\z_2 (w,z) = \exp \sum_{\ell=1}^{\infty} {w^\ell\over \ell} \Xi_\ell(z).
$$
Then we have
$$
\z_2 (1,z) = (1- z)\, \z (\p,z) \quad\hbox{and}\quad \z_2 (w,1) = \z (\psi,w)
$$
wherever the series expansions converge absolutely.
\end{proposition}
We now study the grand partition function $\Xi_\ell(z)$. We first rewrite it in the form
\begin{eqnarray}
\Xi_\ell(z) &=& \sum_{\scriptstyle \s=S^\ell\s }
z^{\sum_{j=0}^{\ell-1}\s_j}\prod_{k=0}^{\ell-1}\psi(S^k\s)\nonumber \\
&=& \sum_{n=0}^{\infty}z^{\ell+n}\sum_{ \s=S^\ell\s \atop \s_0+\cdots +\s_{\ell-1}=\ell+n}\;
 \prod_{k=0}^{\ell-1}\psi(S^k\s) \nonumber
\end{eqnarray}
Notice that for $n=0$ the last sum yields only one term corresponding
to the fixed point $1^{\infty}$ of $S$ (and of $T$).
More generally, the sum over periodic points yields  
${n+\ell -1 \choose \ell -1}={n+\ell -1 \choose n}$ terms, corresponding to the number of
ways of distributing $n$ identical objects into $\ell$ distinct boxes.
\begin{lemma} The sequence $\{ \, \Xi_\ell (z)\, \}$ is uniformly bounded in $\Dc$.
In particular, for each $\ell >0$, 
the power series expansion of $\Xi_\ell (z)$ converges absolutely in $\Dc$.
\end{lemma}
{\it Proof.} Recalling that $\sum_n \Var_n W \leq C_3\t /(1-\t)< \infty$ we have,
for any real $z$ such that $0\leq z \leq 1$,
$$
C^{-1}\, \Lambda_\ell(z) \leq \Xi_\ell (z) \leq C\, 
\Lambda_\ell(z)
$$
where $\Lambda_\ell(z)$ is defined in (\ref{Lambda_n}) and 
$C= e^{{C_3\theta\over 1-\t}}$.
Moreover, from the previous Section it follows that 
${\M}^\ell 1\to h$ in the supremum norm. 
Therefore, since the power series expansion of $\Xi_\ell (z)$ has positive
coefficients, we have for all $z\in \Dc$ and all $\ell >0$, 
$$
|\Xi_\ell (z)| \leq \Xi_\ell (|z|) \leq \Xi_\ell (1) \leq 
C\, \Lambda_\ell(1) \,
\leq  C\, (\Vert h \Vert _{\infty} + o(1)).\qquad \qed
$$
\vskip 0.1cm
\noindent
From the above result we have that for any $z\in \Dc$ the function $\z_2(w,z)$,
viewed as function of the variable $w$, converges absolutely for
$|w|< 1/\exp P(|z|)$. 
To obtain more information we shall establish in the following
theorem a correspondence
between the analytic properties of $\z_2(w,z)$
and the spectral properties of the operator-valued function ${\M}_z$ studied 
in Section \ref{opvalpowser}.
\begin{theorem} The two-variables zeta function
defined in Proposition \ref{twozeta} has the following analytic properties:
\begin{enumerate}
\item for $z\in \Dc$, 
$1/\zeta_2(w,z)$, considered as a
function of the variable $w$, 
is holomorphic in the disk of radius $1/\t \exp P(|z|)$. Its
zeroes in this disk, counted with multiplicity, are the inverses of the
eigenvalues of ${\M}_{z}:{\Ft} \to {\Ft}$ 
in the corresponding
annulus. Moreover, for $0<z\leq 1$, the zero 
of smallest modulus is simple and
located at $1/\exp P(z)$;
\item for $w\in \Dc$, $1/\zeta_2 (w,z)$, considered as a
function of the variable $z$, is holomorphic in the disk of radius
$1$. Its zeroes in this disk are located at 
those values of $z$ such
that ${\M}_{z}:{\Ft}  \to {\Ft}$  has $1/w$ as an eigenvalue.
\end{enumerate}
\end{theorem}
{\sl Proof.} Fixing $n>0$ we let $\sum_{\eta}$ be the sum over
words $\eta$ of length $n$ and denote by
$\s^{(\eta)}$ the periodic concatenation 
$({ \overline {\s_0 \s_1 \dots \s_{n-1} }})$.
Let moreover 
$1_{\eta}\in {\Ft} (\S)$ be such that
$1_{\eta}(\s) = 1$ if $\s$ begins with the word $\eta$, 
$1_{\eta}(\s) = 0$ otherwise. Then we have the following relation for the
grand partition function $\Xi_\ell (z)$:
\be \label{haydn}
\Xi_\ell(z) = \sum_{\scriptstyle \s\in \S \atop \scriptstyle S^\ell\s=\s }
\exp \sum_{j=0}^{\ell-1}W_z(T^j\s) = \sum_{\eta}
({\M}^\ell_{z} 1_{\eta})(\s^{(\eta)}).
\ee
The assertion now follows by
putting together (\ref{haydn}), Theorem \ref{spectrum} and
a straightforward extension of (\cite{Hay}, Theorem 4) to the present situation $\qed$.
\vsni
\noindent
From Proposition \ref{relations} and Proposition \ref{twozeta}
we then have the following,
\begin{corollary} \label{zetas}
\noindent
\begin{enumerate} 
\item $1/\z (\psi,w)$ is holomorphic in the disk
of radius $1/\theta$.
Its zeroes in this disk, counted with multiplicity, are the inverses of
the eigenvalues of ${\M}:\Ft \to \Ft$ in the 
annulus $\{\, \lambda \, : \, \theta < |\lambda | \leq 1\, \}$. 
\item $1/\zeta (\p, z)$ is holomorphic in the disk of radius 
$1$. 
In this disk, $1/\zeta (\p, z)=0$ if and only if $1/z$ is an
eigenvalue of ${\L}: \Ft \to \Ft$.
\end{enumerate}
\end{corollary}
The above result yields 
no zeroes of $1/\zeta (\p, z)$ but the point $z=1$. In this case we know
that the eigenvalue $1$ of ${\L}$ is not isolated
(i.e. there is no `gap'). 
Nevertheless, one may investigate the singular behaviour of $\zeta (\p, z)$ when $z\uparrow 1$.
To this end, consider again eq. (\ref{haydn}) and use (\ref{decomp1}) to rewrite it in the following way:
\begin{eqnarray}
\Xi_\ell(z) 
&=&\l_z^\ell\cdot \sum_{\eta}h_z(\s^{(\eta)})\cdot \nu_z (1_\eta)\, + \, 
\sum_{\eta} ({\Nc}^\ell_{z} 1_{\eta})(\s^{(\eta)})\nonumber \\
&=&\l_z^\ell + \l_z^\ell \cdot \nu_z\left(\sum_{\eta} h_z(\s^{(\eta)}) \cdot 1_{\eta}-h_z\right)
+ \,  \sum_{\eta} ({\Nc}^\ell_{z} 1_{\eta})(\s^{(\eta)})\nonumber \\
&=&:\l_z^\ell +R^{(1)}_\ell (z)+R^{(2)}_\ell (z). \nonumber 
\end{eqnarray}
Now, using the fact that $h_z \in \Ft$ for $z\in J\cup \{1\}$
and reasoning as in the proof of Proposition \ref{specdet2}, one gets
$|R^{(1)}_\ell (z)|\leq C_1 \, \gamma^\ell$ and $|R^{(2)}_\ell (z)|\leq C_2 \, \gamma^\ell$
for $z\in H$ and $\gamma =\max \{ \t, 1-\epsilon\}$ (the notation is as in Section \ref{section4}). 
Therefore putting together (\ref{haydn}), the
above observations, Corollary \ref{singular}, Proposition \ref{twozeta} and Corollary \ref{zetas}, we obtain the following 
(see \cite{Ga}, \cite{Is1} 
for related results):
\begin{theorem} \label{zetina} For all finite $d> -1$
the zeta function $\zeta (\p, z)$ is holomophic in $\D$ and
extends continuously to $\Dc \setminus \{1\}$, whereas in $z=1$ it has 
non-polar singularity. More specifically, it can be written as
$$
\zeta (\p, z) = {A(z)\cdot L(z)\over (1-z)},
$$
where $A(z)$ is as in (\ref{generatingfct1}) and $L(z)$ is continuous on the unit circle $|z|=1$, with $L(1)\neq 0$.
\end{theorem}
\begin{remark} {\rm In the situation considered in Example 1 
a straightforward calculation gives
$$
\Xi_\ell (z) =\left( \sum_{n=1}^{\infty}z^{n}p_n\right)^\ell\quad\hbox{and}\quad
\z_2 (w,z) = \left(1- w\sum_{n=1}^{\infty}z^{n}p_n\right)^{-1}.
$$
Using Proposition \ref{twozeta} we then get for $\zeta (\p, z)$ the expression given in Theorem \ref{zetina} with 
$A(z)=\left(1- \sum_{n=1}^{\infty}z^{n}p_n\right)^{-1}$ and $L(z)\equiv 1$.}
\end{remark}
\begin{remark} {\rm A singularity of non-polar nature in the $\zeta$-function should perhaps be
related to the presence of some kind of `phase transition' in the system (see remark \ref{weakgibbs} above and also \cite{PS}).
This is the case
for instance in \cite{Ga} where a direct connection with the model introduced by Fisher \cite{FF} is emphasized.}
\end{remark}
\vsni 
\section{Appendix: interval maps with indifferent fixed points}\label{app}
We shall consider a class of non-uniformly expanding interval maps $F : \ui \to \ui$ satisfying
the following assumptions:
\begin{enumerate}

\item there is a number $q\in (0,1)$ such that the restrictions $F|_{(0,q)}$ and $F|_{(q,1)}$ 
extend to $C^1$-diffeomorphisms on $I_0=[0,q]$ and $I_1=[q,1]$, which are $C^2$ for $x>0$, and such that
$F(0)=0$ and $F(I_0)=F(I_1)=\ui$;

\item there are numbers $0 < \beta < 1$ and $ m \in \Z^+$
such that $(F^\ell)' \geq 1/\beta$ on $I_1$ for all $\ell \geq m$;
whereas $F'> 1$
on $(0,q)$ and $F' (0) = 1$; 

\item $F$ has the following asymptotic behaviour when $x\to 0_+$:
$$
F(x) =x + r\, x^{1+s}(1+u(x))
$$
for some constant $r >0$, exponent $1+s>1$ and where $u(x)$ satisfies
$u(0)=0$ and $u'(x)={\cal O}(x^{t-1})$ as $x\to 0_+$
for some $t>0$.
\end{enumerate}
The partition of $\ui$ whose elements are the intervals 
$I_0$ and
$I_1$ is a Markov partition for $F$. Let $\O$ be as in Section \ref{inducing}. The map $\pi:\O \rightarrow \ui$
defined by 
\be
\pi(\o) =x\quad\hbox{according to}\quad F^j(x)\in I_{\o_j},\;\; j\geq 0
\ee  
is a coding map which is a homeomorphism on the residual set of points in $[0,1]$ 
which are not preimages of $1$ with the map $F$. 
Moreover $F\circ \pi = \pi \circ T$. Let $F_i$ be
the inverse branch of $F$ on $I_i$, $i=0,1$.  Given $\o \in \O_0$ 
we set
\be\label{pot}
V(\o) =\log [(F_{\o_0})^{\prime}(\pi(\o_1 \o_2\dots ))]  .
\ee
It is then easy to realize that its induced version $W(\o)$ can be written as
\be\label{indpot}
W(\o) = \log [(G_{\s_0(\o)})^{\prime}(\pi (\s_1(\o) \s_2(\o)\dots ))],
\ee
where $\s_0(\o) \s_1(\o)\dots =\iota ( \o)$ (see Section \ref{opvalpowser}) and
$G_i$, $i\geq 1$, are the (countably many) inverse branches of the induced map
\be
x\rightarrow G(x) = F^{\tau (x)}(x),
\ee
with
\be\label{fpf1}
\tau(x)=1+\min \{n\geq 0 \;:\; F^n(x)\in A_1\;\}.
\ee
For notational simplicity' sake we shall denote with the same symbol the first passage function (\ref{fpf}) 
and its lift (\ref{fpf1})
with $\pi$, as well as the level sets
\be
A_n=\{x\in \ui : \tau(x) = n\} = [F_0^{n}(1),F_0^{n-1}(1)]\, .
\ee
Notice that $G$ (once suitably extended to the set $\{F_0^{n}(1)\}_{n\geq 0}$) maps $A_n$
onto $\ui$, for all $n\geq 1$. 

\noindent
It turns out that the overall statistical behaviour of the map $F$ depends 
on the way the lengths $|A_n|$ of the levelsets $A_n$ vanish as $n\to \infty$.
\begin{lemma}\label{cienne}
Under the hypotheses (1)-(3) on the map $F$ we have for $n\to \infty$
$$
|A_n| \sim  r(rsn)^{-1-1/s}\; .
$$
\end{lemma}
{\it Proof.} The last property of $T$ gives for the inverse
function:
$$  
F_0(x) = x - r x^{1+s}(1+v(x))
$$ 
where $v(x)$ is such that $v(0)=0$ and $v'(x)={\cal O}(x^{t-1})$ as $x\to 0_+$.
We write this expression in a
more manageable form, that is:
$$   
F_0(x) =\biggl( x^{-s} + rs (1+\vt(x)) \biggr)^{-1/s}
$$ 
where $\vt (x)$ is another function such that $\vt (0)=0$ and
$v^{\prime}(x)-\vt^{\prime}(x) ={\cal O}(x^{u-1})$ where $u=\min \{s,t\}$. It is then easy
to check that
$$
F_0^n(x) =  \left(x^{-s} + r sn\biggl(1+
{1\over n}\sum_{l=0}^{n-1}\vt(F_0^l(x)) 
\biggr)\right)^{-1/s}
$$
whence
\be\label{cala}
F_0^n(1) =  (rs n)^{-1/s}\left(1+ o(1))\right)^{-1/s}
\ee
and the assertion follows. $\qed$

\noindent
We now list some properties of the induced map $G$ which are relevant for our discussion.
\begin{proposition}\label{G}\
\begin{enumerate}
\item $G_{|_{A_n}}$ is a $C^2$-diffeomorphism of $A_n$ onto $\ui$,
for all $n\geq 1$;
\item $\exists m \in \Z^+$ so that
$$ 
\inf_{\scriptstyle x\in \A_n \atop \scriptstyle n\geq 1}|(G^m)'(x)| = 1/\beta > 1;
$$
\item
$$
\inf_{\scriptstyle x,y,z\in \A_n \atop \scriptstyle n\geq 1} 
\left| {G''(x)\over G'(y)G'(z)}\right| = K < \infty .
$$
\end{enumerate}
\end{proposition}
{\it Proof.} Statements 1) and 2) are immediate consequences of the
definition.
To show 3) we first observe that the chain rule yields
$$
{G''(x)\over (G')^2(x) } =
\sum_{k=0}^{\tau (x)-1}{F''(F^k(x))\over (F')^2(F^k(x))}\cdot 
{1\over \prod_{j=k+1}^{\tau (x)-1}F'(F^j(x))}.
$$
On the other hand, the properties of $F$ imply that
$$
{|F''(x)|\over |(F')^2(x)|}  \asymp  x^{s-1}
$$
Moreover, if $\xi_n$ is any point in 
$A_{n}$, then $\xi_n^{s-1} \asymp n^{-1+{1\over s}}$ and 
$\prod_{j=0}^{n-1}|F'(F^j(\xi_n))| \equiv |G'(\xi_n)| \asymp n^{1+{1\over s}}$.
Putting together these remarks we get
$$
\left|{G''(\xi_n)\over (G')^2(\xi_n) }\right| \leq C_1 \, \sum_{k=0}^{n-1} (n-k)^{-2}
\leq C_2  \, \sum_{k=1}^{\infty} k^{-2} \leq C_3,
$$
and the assertion follows by noting that $G'(\xi_n) \asymp G'(\eta_n)$ for any choice of
$\xi_n,\eta_n\in A_n$ and any $n\geq 1$. $\qed$
\vskip 0.1cm
\noindent
The above properties yield a uniform bound for the buildup of non-linearity
in the induction process. 
\begin{corollary} Let $x,y \in \ui$ be such that
$G^j(x)$ and $G^j(y)$ belong to the same partition set $A_{k_{j}}$, for $0\leq j \leq n$
and some $n\geq 1$. 
Then there are two constants $C>0$ and $\alpha <1$ such that
$$
\left|\log {G'(x)\over G'(y)}\right| \leq C\, \alpha^n\, .
$$
\end{corollary}
{\it Proof.} 
Taking $x,y\in A_{k_0}$, let $\eta\in A_{k_0}$ be such 
that $|G'(\eta)|=|A_{k_0}|^{-1}$. Then, using Proposition \ref{G}(3), we have
\begin{eqnarray}
\left| \, \log {G'(x)\over G'(y)}\, \right| 
&=&\left|{G''(\xi)\over G'(\xi) }\right|\cdot |x-y|
\quad\hbox{for some}\quad \xi \in [x,y] \subseteq A_{k_0} \nonumber \\
&=& \left|{G''(\xi)\over G'(\xi)G'(\eta) }\right|\cdot {|x-y|\over |A_{k_0}| } 
\leq K\, {|x-y|\over |A_{k_0}| } \, \cdot\nonumber
\end{eqnarray}
Now, using Proposition \ref{G}(2),
we can find a constant $C_4>0$ such that, under the above hypotheses,
$|x-y| \leq C_4\, |A_{k_0}| \, \alpha^n$ with $\alpha = \beta^{1\over m}<1$.
$\qed$

\noindent
These results and (\ref{indpot}) imply that the potential $V$ defined in (\ref{pot}) satisfies the properties
listed in Section \ref{inducing} for every $\theta \geq \alpha$. We can then apply the whole subsequent theory. 
In particular, there is an unique absolutely continuous 
probability measure $\rho (dx) =  h(x)\, dx$ which is invariant
for the dynamical system $(\ui, G)$ and whose density 
$h$ is Lipschitz continuous and satisfies $h \asymp 1$ (see also \cite{Wal2}). In turn, this and Lemma \ref{cienne} imply 
(see, e.g., \cite{CI2}, Lemma 2.4)
that as $n\to \infty$,
\be\label{levset}
\rho (A_n) \sim \,h(0)\,  |A_n|\sim h(0)\, r(rsn)^{-1-1/s}\; .
\ee
Moreover, the $\s$-finite absolutely continuous measure 
$\mu (dx) = e(x)\, dx$ with
\be
e =\sum_{n=0}^{\infty}h\circ F_0^{n}\cdot (F_0^{n})^{\prime}
\ee
is invariant for $(\ui , F)$.  
It can be shown (see \cite{Th}) that $e(x) \asymp x^{-s}$. Clearly, both $e$ and $h$ are the
lifts with the map $\pi$ of the corresponding quantities considered e.g. in Section \ref{section4}.
From (\ref{levset}) and Definition \ref{degree} it follows that $([0,1],F,\mu)$ has ergodic degree 
\be
d={1\over s}-1.
\ee
However, the asymptotic equivalence expressed by (\ref{levset}) is somewhat stronger than the mere knowledge
of the ergodic degree. Indeed, (\ref{levset}) implies that all symbols $\approx$ in Proposition
\ref{generfcts} can be turned into true asymptotic equivalences $\sim$. In particular, application of a Tauberian
theorem for power series to the functions $A(z)$ and $B(z)$ is now more informative and gives
\be
a_n \sim \cases{ c_1\, n^{-1+1/s} &if $s>1$, \cr
                 c_2\, (\log n)^{-1} &if $s=1$, \cr }
\ee
and 
\be 
b_n \sim c_3 \, n^{1-1/s}\;\;\hbox{if}\;\; 0<s<1 ,
\ee
for some positive constants $c_1,c_2,c_3$. According to Theorem \ref{zetina} the behaviour of $a_n$ given above 
can be used to establish (again by a Tauberian theorem for power series) the precise asymptotic 
behaviour when $z\uparrow 1$ of the dynamical zeta function for the map $F$.
However we shall leave this easy task to the interested reader and end this Section by discussing
some consequences of the construction outlined in the main part of the paper on scaling and mixing properties of
$(\ui, F, \mu)$. To this end we set
\be
 B_{+}:= \cup_{\epsilon} \{E\in {\cal B}(\ui): \, m (E)>0, 
\, E \subseteq \ui \setminus  (0, \epsilon)\,\}\, .
\ee
The following sharpening of Theorem \ref{mr} is then a straightforward consequence of application of the arguments 
of Section \ref{rensca}
to this situation. 
\begin{theorem}\label{poly} For all $E\subset B_+$ we have
\begin{itemize}
\item ${\mu (E\cap T^{-n}E)/ (\mu(E))^2} \sim c_1\, n^{-1+1/s}$ if $s> 1$;
\item ${\mu (E\cap T^{-n}E)/ (\mu(E))^2} \sim c_2\, (\log n)^{-1}$ if $s=1$;
\item $[{\hat \mu} (E\cap T^{-n}E)-({\hat \mu}(E))^2]/({\hat \mu}(E))^2  \sim c_3 \, n^{1-1/s}$ if $0<s<1$;
\end{itemize}
where, for $s<1$, we have set ${\hat \mu}=\mu/\mu(\ui)$.
\end{theorem}
We finally concentrate on the finite measure case $s<1$. Given a function $f:\ui \to \R$ which is H\"older continuous 
with exponent $\gamma >0$, that is:
$|f(x)-f(y)|\leq K |x-y|^\gamma$, let us denote by ${\tilde f}$ its projection $f\circ \pi:\O \to \R$.
It is easy to realize that, for $n$ large, 
\be 
\var_n {\tilde f} = \sup_{x \in A_{n+1}}|f(x)-f(0)|\leq 2K (rs n)^{-\gamma/s}
\ee
where the last inequality comes from the H\"older property along with (\ref{cala}). Therefore,
in order to have ${\tilde f}\in {\cal H}_a$ with $a\geq d$ it is sufficient that $\gamma \geq 1-s$. Thus,
using Theorem \ref{pcp} and
taking into account that due to Theorem \ref{poly} the bounds for the weak-Bernoulli property are truly
polynomial (no slowly varying functions being involved) we have the following (see \cite{Yo}, \cite{LSV1},
\cite{Hu}, \cite{Sa2} for related
results):
\begin{theorem}\label{clus} Let $0<s<1$ and $f,g:\ui \to \R$ be H\"older continuous with exponent $\gamma \geq 1-s$. Then
$$
|{\hat \mu} (f\cdot g\circ F^{n}) - {\hat \mu}(f){\hat \mu}(g)|\, = \, {\cal O} (n^{1-1/s})\, .
$$
\end{theorem}


\begin{thebibliography}{DEGHL}


\bibitem[Aa]{Aa}
{\sc J Aaronson}, {\it An introduction to infinite ergodic theory}, AMS, 1997.

\bibitem[Ab]{Ab}
{\sc L M Abramov}, \,  {\sl The entropy of a derived automorphism}, 
Amer. Math. Soc. Transl. (2) {\bf 49} (1965), 162-166. 

\bibitem[Ba1]{Ba1}
{\sc V Baladi}, \, {\sl Positive Transfer Operators and Decay of Correlations}, World Scientific, 2000. 

\bibitem[Ba2]{Ba2}
{\sc V Baladi}, \, {\sl Dynamical zeta functions}, Real and Complex
Dynamical Systems (B. Branner and P. Hjorth eds.), Kluwer Academic
Publishers, 1995. 

\bibitem[Bo]{Bo}
{\sc R Bowen}, \, {\it Equilibrium states and the ergodic theory of
Anosov diffeomorphisms},  LNM 470 (1975), Springer-Verlag.

\bibitem[Br]{Br}
{\sc X Bressaud}, \,  {\sl Subshift on an infinite alphabet}, Erg. Th. Dyn. Syst. {\bf 19} (1999), 1175-1200. 

\bibitem[CI]{CI}
{\sc M Campanino, S Isola}, \,  {\sl Statistical properties of long return
times in type I intermittency}, Forum Mathem. {\bf 7} (1995), 331-348. 

\bibitem[CI2]{CI2}
{\sc M Campanino, S Isola}, \,  {\sl Infinite invariant measures or non-uniformly expanding transformations of $\ui$: weak law
of large numbers with anomalous scaling}, Forum Mathem. {\bf 8} (1996), 71-92. 

\bibitem[Che]{Che} {\sc N Chernov}, {\sl Limit theorems and markov approximations for
chaotic dynamical systems}, Probab. Theory Relat. Fields {\bf 101},
(1995) 321-362.

\bibitem[Chu]{Chu}
{\sc K L Chung}, {\it Markov chains with stationary
transition probabilities}, Springer 1967.

\bibitem[Fe]{Fe}
{\sc W Feller}, \,  {\sl An Introduction to Probability Theory and Its
Applications}, Volume 2,
J.Wiley and Sons, New York 1970.

\bibitem[FF]{FF}
{\sc B U Federhof, M E Fisher},
Annals of Physics (N.Y.) {\bf 58} (1970).

\bibitem[FL]{FL}
{\sc A M Fisher, A Lopes},
{\sl Exact bounds for the polynomial decay of correlations, 
$1/f$ noise and the CLT for the equlibrium state of a non-H\"older potential},
Nonlinearity  {\bf 14} (2001), 1071-1104. 

\bibitem[Ga]{Ga}
{\sc G Gallavotti}, \, {\sl Funzioni zeta e insiemi basilari}, 
Accad. Lincei Rend. Sc. fis. mat. e nat. {\bf 61} (1976), 309-317.

\bibitem[Har]{Har}
{\sc G H Hardy}, {\it Divergent series},  Oxford at the Calrendon Press 1949.

\bibitem[Hay]{Hay} {\sc N T A Haydn},
{\sl Meromorphic extension of the zeta
function for Axiom A flows},  Erg. Th. Dyn. Sys. {\bf
10} (1990), 347-360.

\bibitem[Ho]{Ho} {\sc F Hofbauer},
{\sl Examples for the nonuniqueness of the equilibrium state},  Trans. Amer. Math. Soc. {\bf
228} (1977), 223-241. 


\bibitem[Hu]{Hu} {\sc H Hu},
{\sl Decay of correlations for maps with indifferent fixed points},  Preprint (1999).

\bibitem[HI]{HI} {\sc N T A Haydn, S Isola},
{\sl Parabolic rational maps},  J. London Math. Soc. (2) {\bf 63} (2001), 673-689.

\bibitem[Is1]{Is1} {\sc S Isola},
{\sl Renewal sequences and intermittency}, 
 J. Stat. Phys. {\bf 97} (1999), 263-280.

\bibitem[Is2]{Is2} {\sc S Isola},
{\sl On the rate of convergence to equilibrium for countable ergodic
Markov chains}, to appear in Markov Proc. Rel. Fields (2002).

\bibitem[Ka]{Ka} {\sc T Kato},
{\it Perturbation theory of linear operators}, 
Springer-Verlag, Berlin Heidelberg New York (1980).

\bibitem[Ke]{Ke} {\sc G Keller}, {\sl On the rate of convergence to equilibrium 
in one-dimensional systems}, 
Comm. Math. Phys. {\bf 96} (1984), 181-193.

\bibitem[Ki]{Ki} {\sc J F C Kingman}, {\it Regenerative Phenomena}, 
John Wiley, 1972.

\bibitem[LSV1]{LSV1}
{\sc C Liverani, B Saussol, S Vaienti},
{\sl Conformal measure and decay of correlations for covering weighted systems},
Erg. Th. Dyn. Sys. {\bf 18} (1998), 1399-1420.


\bibitem[LSV2]{LSV2}
{\sc C Liverani, B Saussol, S Vaienti},
{\sl A probabilistic approach to intermittency},
Erg. Th. Dyn. Sys. {\bf 19} (1999), 671-685.

\bibitem[MD]{MD}
{\sc V Maume-Deschamps},
{\sl Correlation decay for Markov maps on a countable state space},
Erg. Th. Dyn. Sys. {\bf 21} (2001), 165-196.


\bibitem[MRTVV]{MRTVV}
{\sc C Maes, F Redig, F Takens, A Van Moffaert, E Verbitsky},
{\sl Intermittency and weak Gibbs states},
Nonlinearity  {\bf 13} (2000), 1681-1698. 


\bibitem[Nu]{Nu} {\sc R Nussbaum}, {\sl The radius of the essential spectrum}, Duke Math.
J. {\bf 37} (1970), 473-478.  

\bibitem[Pol1]{Pol1} {\sc M Pollicott}, {\sl Meromorphic extensions of generalised
zeta functions},  Invent. math. {\bf 85} (1986), 147-164.

\bibitem[Pol2]{Pol2} {\sc M Pollicott}, {\sl Rates of mixing for potentials of summable variation}, 
Trans. Amer. Math. Soc. {\bf 352} (2000), 843-853.

\bibitem[PP]{PP}
{\sc W Parry, M Pollicott}, \, {\it Zeta functions and the
periodic orbit structure of hyperbolic dynamics}, Soci\'et\'e
Math\'ematique de France (Ast\'erisque {\bf 187-188}), Paris. 

\bibitem[Pos]{Pos}
{\sc A G Postnikov}, \, {\it Tauberian Theory and its Applications}, Proceedings of the Steklov Institute
of Mathematics, 1980, Issue 2. 

\bibitem[PS]{PS} {\sc T Prellberg, J Slawny}, {\sl Maps of intervals with 
indifferent fixed points: thermodynamic formalism and phase transitions},
J. Stat. Phys. {\bf 66} (1992), 503-514.

\bibitem[RS]{RS}
{\sc M Reed, B Simon}, \,{\it Methods of Modern Mathematical Physics}, Vol. IV:
Analysis of Operators, Academic Press, New York 1978.

\bibitem[Ru1]{Ru1}
{\sc D Ruelle}, \,{\sl Zeta functions for expanding maps and
Anosov flows}, Invent. Math. {\bf 34} (1976), 231-242.

\bibitem[Ru2]{Ru2}
{\sc D Ruelle}, {\sl Dynamical Zeta Functions for Piecewise
Monotone Maps of the Interval}, 
American Mathematical Society (CRM Monograph Series, {\bf 4}), 
Providence, Rhode Island USA, 1994.

\bibitem[Ru3]{Ru3}
{\sc D Ruelle}, \,{\sl One dimensional Gibbs' states 
and Axiom A diffeomorphisms}, J. Diff. Geom. 
{\bf 25} (1987), 117-137.

\bibitem[Ru4]{Ru4}
{\sc D Ruelle}, \,{\it Thermodynamic Formalism}, 
 Addison-Wesley Publ. Co. 1978.

\bibitem[Sa1]{Sa1}
{\sc O Sarig},
{\sl Phase transitions for countable Markov shifts},
Commun. MAth. Phys. {\bf 217} (2001), 555-577.

\bibitem[Sa2]{Sa2}
{\sc O Sarig},
{\sl Subexponential decay of correlations},
Preprint 2001.

\bibitem[Si]{Si}
{\sc B Simon}, \,{\it The Statistical Mechanics of Lattice Gases}, 
Princeton University Press, 1993.

\bibitem[Th]{Th} {\sc M Thaler }, {\it Estimates of the invariant densities of endomorphisms
with indifferent fixed points}, Israel Jour. Math. {\bf 37} (1980), 303-314.


\bibitem[Wal1]{Wal1} {\sc P Walters}, {\sl Ruelle's operator theorem and $g$-measures}, Trans. Amer. Math. Soc. 
{\bf 214} (1975), 375-387.
 

\bibitem[Wal2]{Wal2} {\sc P Walters}, {\sl Invariant measures and equilibirum states for some mappings
which expand distances}, Trans. Amer. Math. Soc. {\bf 236} (1978), 121-153.
 

\bibitem[Yo]{Yo} {\sc L S Young}, {\sl Recurrence times and rate of mixing},
 Isr. J. Math. {\bf 110} (1999),  153-188. 

\bibitem[Yu]{Yu} {\sc M Yuri}, {\sl Thermodynamic formalism for certain non-hyperbolic maps},
Erg. Th. Dyn. Sys. {\bf 19} (1999), 1365-1378. 


\bibitem[Zig]{Zig}
{\sc A Zigmund}, \,{\it Trigonometric Series}, Cambridge at the University Press, 1968.

\end{thebibliography}
\end{document}